\documentclass[11pt,a4paper]{article}
\usepackage[T1]{fontenc}
\usepackage[utf8]{inputenc}
\usepackage{lmodern,microtype}
\usepackage[a4paper,margin=0.92in]{geometry}
\usepackage{amsmath,amssymb,amsthm,mathtools,mathrsfs,bm}
\usepackage{enumitem,booktabs,longtable,array}
\usepackage{xcolor,hyperref}
\usepackage[nameinlink,noabbrev]{cleveref}
\hypersetup{colorlinks=true,linkcolor=blue!55!black,citecolor=blue!55!black,urlcolor=blue!55!black,
 pdftitle={Min--Max Theory for Surfaces with Prescribed Anisotropic Mean Curvature},pdfauthor={Anonymous manuscript}}
\allowdisplaybreaks
\numberwithin{equation}{section}
\setlist[enumerate]{leftmargin=2em,itemsep=2pt,topsep=4pt}
\newtheorem{theorem}{Theorem}[section]
\newtheorem{proposition}[theorem]{Proposition}
\newtheorem{lemma}[theorem]{Lemma}
\newtheorem{corollary}[theorem]{Corollary}

\theoremstyle{definition}
\newtheorem{definition}[theorem]{Definition}

\newcommand{\R}{\mathbb R}
\newcommand{\dd}{\,\mathrm d}\newcommand{\eps}{\varepsilon}
\newcommand{\Vol}{\operatorname{Vol}}
\newcommand{\dist}{\operatorname{dist}}\newcommand{\spt}{\operatorname{spt}}

\newcommand{\an}{\operatorname{An}}\newcommand{\tr}{\operatorname{tr}}

\newcommand{\pair}[2]{\left\langle#1,#2\right\rangle}
\newcommand{\calF}{\mathcal F}\newcommand{\bfF}{\mathbf F}
\newcommand{\Fphi}{\mathbf F_{\Phi}}
\newcommand{\Aph}{\mathcal A_{\Phi,h}}\newcommand{\Wphi}{\mathbf W_{\Phi}}
\newcommand{\Wph}{\mathbf W_{\Phi,h}}
\newcommand{\Cacc}{\mathcal C(M)}\newcommand{\tPhi}{\widetilde\Phi}

\newcommand{\one}{\mathbf 1}
\title{Min--Max Theory for Surfaces with Prescribed Anisotropic Mean Curvature}
\author{XIAOXIANG JIAO AND WENDUO ZOU}
\date{}

\begin{document}
\maketitle
\begin{abstract}
Let $(M^3,g)$ be a closed Riemannian manifold, and let
$\Phi$ be a smooth even elliptic integrand. In this paper,
we develop a one-parameter min--max theory for the functional
$\mathcal A_{\Phi,h}(E)
=\mathbf F_\Phi(E)-\int_E h\,\dd V_g$,
where $h\in C^\infty(M)$. As an application, we establish
the existence of a nonempty, smooth, closed, almost embedded
surface of prescribed anisotropic mean curvature $h$,
arising as the boundary of a Caccioppoli set, for an open
and dense subset of $C^\infty(M)$.

Our proof builds on the min--max techniques developed in
\cite{DD24,DDL24,ZZ20} and adapts them to the prescribed
anisotropic setting. A key new ingredient is a jet
transversality argument that identifies an open and dense
class of prescription functions for which the min--max
regularity theory can be established. 
\end{abstract}
\medskip

\section{Introduction}
\label{sec:introduction}

Let $(M^3,g)$ be a closed Riemannian manifold and let $\Phi$
be a smooth elliptic integrand on $M$. The anisotropic energy
of a smooth surface $\Sigma$ is defined by
\[
\mathbf F_\Phi(\Sigma)
=\int_\Sigma\Phi(x,\nu)\,\dd\mathcal H_g^2.
\]
Here $\Phi$ is a positive even function on the unit tangent bundle
$UM$ whose positively one-homogeneous extension is elliptic.
Critical points of this energy are called anisotropic
minimal surfaces. More generally, given a smooth function
$h\in C^\infty(M)$, we consider surfaces with prescribed anisotropic
mean curvature $h$. Such surfaces arise as critical boundaries
of the functional
\[
\mathcal A_{\Phi,h}(E)
=\mathbf F_\Phi(E)-\int_E h\,\dd V_g,
\]
defined on the space of Caccioppoli sets in $M$. With our sign
convention, their anisotropic mean curvature is $h$ with respect
to the outward unit normal of the set.

The existence of closed minimal hypersurfaces in Riemannian
manifolds is a fundamental problem in geometric analysis.
The min--max theory initiated by Almgren and developed by
Pitts \cite{Pit81}, together with the regularity theory of
Schoen and Simon \cite{SS81}, gives the existence of a nontrivial
smooth closed embedded minimal hypersurface in every closed
Riemannian manifold of dimension between three and seven.
In recent years, this theory has been further developed in
several directions. In particular, Marques and Neves \cite{MN16}
established Morse index bounds for min--max minimal hypersurfaces
and formulated the multiplicity one conjecture as part of their
program to develop a Morse theory for the area functional.

An important extension of the classical min--max theory was
given by Zhou and Zhu \cite{ZZ19}, who developed a min--max
theory for hypersurfaces with constant mean curvature.
They subsequently extended this theory to hypersurfaces with
prescribed mean curvature and proved existence for a generic
class of smooth prescription functions \cite{ZZ20}.
Besides providing new existence results, the prescribed mean
curvature theory plays an important role in the study of minimal
hypersurfaces. By approximating the minimal hypersurface problem
by prescribed mean curvature problems, Zhou \cite{Zho20} proved
the multiplicity one conjecture for bumpy metrics. Together
with the work of Marques and Neves \cite{MN21}, this led to
precise Morse index results for minimal hypersurfaces realizing
the volume spectrum for generic metrics.

The existence and regularity theory for critical points of
anisotropic surface energies has a long history in geometric
measure theory. While the theory for minimizers of elliptic
integrands is well developed, the study of general critical
points is substantially more difficult. One of the main
difficulties is the lack of a general monotonicity formula,
which plays a central role in the regularity theory of
isotropic minimal surfaces. Important early contributions
to the anisotropic theory include the regularity results of
Allard \cite{All83,All86} and the curvature and compactness
results of White \cite{Whi87a}. More recently, De Philippis,
De Rosa and Ghiraldin \cite{DDG18} extended Allard's
rectifiability theorem to varifolds with locally bounded
anisotropic first variation under the atomic condition.
De Rosa and Tione \cite{DT22} obtained almost-everywhere
regularity for Lipschitz graphs with anisotropic mean curvature
in $L^p$, for $p$ greater than the dimension of the graph,
under a suitable ellipticity condition.

A major step in the development of anisotropic min--max
theory was made by De Philippis and De Rosa \cite{DD24},
who established the existence of anisotropic minimal surfaces
and surfaces with constant anisotropic mean curvature in
closed Riemannian three-manifolds, with at most one singular
point. Subsequently, De Philippis, De Rosa and Li \cite{DDL24}
extended the existence theory to arbitrary dimensions (at least
three) and obtained improved regularity,thereby answering a conjecture of Allard
\cite[p.~288]{All83}. In particular, they obtained smooth closed
anisotropic minimal surfaces in dimension three. The key idea
in their work is to construct optimal nested volume-parameterized
sweepouts with uniform mass-ratio bounds, which provide the local
area control needed in the regularity theory.

Several other approaches to anisotropic min--max theory
have recently been developed. De Rosa and Pigati \cite{DRP25}
introduced a phase-transition approach based on the anisotropic
Allen--Cahn energy, producing nontrivial integral varifolds
stationary with respect to anisotropic surface energies.
De Rosa, Halavati and Wang \cite{DHW26} established an anisotropic
version of the Meeks--Simon--Yau theorem and developed an
anisotropic Simon--Smith min--max theory with genus bounds.
They also proved removable singularity theorems under suitable
assumptions on the integrand. De Rosa and Halavati \cite{DH26}
established a Weyl law for the codimension-one anisotropic
volume spectrum and obtained equidistribution in the
Grassmannian bundle by varying the anisotropy while keeping
the ambient metric fixed. These results show that several
important aspects of the classical min--max theory are now
being extended to the anisotropic setting.

There has also been recent progress in the regularity theory
for varifolds with bounded anisotropic first variation.
De Philippis and Pigati \cite{DPP26} established Michael--Simon
inequalities for classes of anisotropic integrands, including
surface energies close to the area functional in dimension
three, and studied their relation to varifold compactness.
Firester and Tsiamis \cite{FT26} extended these inequalities
to a much broader class of integrands, including convex even
hypersurface anisotropies. De Rosa, Firester and Tsiamis
\cite{DFT26} subsequently established an anisotropic Allard-type
regularity theorem in arbitrary codimension under suitable
assumptions on the integrand. For a recent overview of the
existence and regularity theory of anisotropic minimal surfaces,
we refer the reader to \cite{DR26}.

In this paper, we extend the anisotropic min--max theory to
surfaces with prescribed anisotropic mean curvature in dimension
three, thereby addressing a natural existence problem in the
calculus of variations. Beyond its own geometric interest, this
problem is closely related to the study of anisotropic minimal
surfaces. As discussed above, in the isotropic setting the
prescribed mean curvature min--max theory was a key ingredient
in the proof of the multiplicity one conjecture. The theory
developed in this paper provides a foundation for the study
of multiplicity one and Morse index in the anisotropic setting.

\begin{samepage}
Our main result is the following.

\begin{theorem}[Existence of anisotropic PMC surfaces]
\label{thm:apmc-existence}
Let $(M^3,g)$ be a closed Riemannian manifold
and let $\Phi$ be a smooth even elliptic integrand on $M$.
There exists an open and dense subset
$\mathcal G_\Phi\subset C^\infty(M)$
such that for every $h\in\mathcal G_\Phi$, there exists a nontrivial,
smooth, closed, almost embedded surface $\Sigma=\partial\Omega$
with prescribed anisotropic mean curvature $h$ with respect to
the outward unit normal of the Caccioppoli set $\Omega$.
\end{theorem}
\end{samepage}

An important step in establishing a min--max theory for
prescribed anisotropic mean curvature is to identify a suitable
class of prescription functions for which the regularity theory
can be carried out. In the isotropic setting, Zhou and Zhu
\cite{ZZ20} introduced a class of Morse functions with
regular zero level sets whose mean curvature vanishes only to
finite order. They proved the density of this class using mean
curvature flow and unique continuation for the mean curvature
evolution equation.

Here we take a different approach. For a fixed elliptic integrand
$\Phi$, we introduce a class of functions, called $\Phi$-good
functions. In addition to the Morse condition and the regularity
of the zero level set, we require that zero be a regular value
of the anisotropic mean curvature function on that level set.
By expressing the anisotropic mean curvature of a regular level
set in terms of the two-jet of its defining function, we use
the jet transversality theorem to prove density. We show that
the set of $\Phi$-good functions is $C^3$-open and $C^\infty$-dense
in $C^\infty(M)$ (Proposition~\ref{thm:good-generic}). This condition
allows us to exclude infinite-order contact between distinct
touching sheets and to establish the touching-set estimates
needed in the regularity theory.

Another main step is to obtain uniform mass-ratio bounds in the
min--max construction. Following \cite{DDL24}, we extend the
construction of optimal nested volume-parameterized sweepouts
to the prescribed functional $\mathcal A_{\Phi,h}$
(Theorem~\ref{thm:MR-ONVP}). With these bounds, we establish the
existence of almost-minimizing varifolds and develop the
corresponding replacement theory. Combining the touching-set
estimates, stable curvature estimates and compactness results
with the replacement and gluing arguments, we obtain the
regularity of the min--max surfaces and prove
Theorem~\ref{thm:apmc-minmax}. Theorem~\ref{thm:apmc-existence}
then follows from Corollary~\ref{cor:APMC-minmax} and the
genericity of $\Phi$-good functions.

The paper is organized as follows.
In Section~\ref{sec:preliminaries}, we introduce the notation
and recall the basic facts concerning anisotropic energies
and their variations.
In Section~\ref{sec:good-local}, we establish the genericity
of $\Phi$-good functions and study the touching sets of
almost embedded surfaces with prescribed anisotropic mean
curvature.
In Section~\ref{sec:transfer}, we develop the min--max
construction and prove the main results.
The appendices contain the derivation of the anisotropic
mean curvature formula and a local regularity result for
anisotropic varifolds under a uniform lower mass-ratio bound.

\vspace{10pt}
\medskip
\noindent\textbf{Acknowledgements.}
Supported by the National Natural Science Foundation of China
(Grant no.~12371055).

\section{Geometric and variational preliminaries}
\label{sec:preliminaries}

Throughout this paper, we fix a smooth closed three-dimensional
Riemannian manifold $(M^3,g)$. In this section, we introduce the
notation and recall the basic definitions and variational
properties used in the paper.

\subsection{Notation and terminology}

For background on geometric measure theory, we refer the reader
to  \cite{Sim83}. We use the following notation:
\begin{itemize}[
    leftmargin=1.2em,
    itemsep=2pt,
    parsep=0pt,
    topsep=3pt,
    partopsep=0pt
]

\item $B_r(p)$ denotes the open geodesic ball of radius $r$
centered at $p$.

\item $\operatorname{inj}(M)$ denotes the injectivity radius
of $(M,g)$.

\item $\Vol$, $\dd V_g$, and $\mathcal H_g^k$ denote the
Riemannian volume, the Riemannian volume measure, and the
$k$-dimensional Hausdorff measure, respectively.

\item $\operatorname{spt}\mu$ denotes the support of a measure
$\mu$, and $\mu\llcorner E$ denotes its restriction to a
measurable set $E$.

\item $U\Subset V$ denotes that $\overline U$ is compact
and contained in $V$.

\item $E\triangle F$ denotes the symmetric difference of
the sets $E$ and $F$.

\item $\mathfrak X(M)$ denotes the space of smooth vector
fields on $M$.

\item $UM=UT(M)$ denotes the unit tangent bundle:
\[
UM:=\{(x,\nu)\in TM:|\nu|_g=1\}.
\]

\item $G_2(M)$ denotes the unoriented two-plane bundle:
\[
G_2(M):=\{(x,T):x\in M,\ T\subset T_xM,\ \dim T=2\}.
\]

\item $\Cacc$ denotes the space of Caccioppoli sets in $M$,
identified up to $\mathcal H_g^3$-null sets and endowed
with the $L^1$ topology:
\[
\|\one_E-\one_F\|_{L^1(M)}=\Vol(E\triangle F).
\]

\item $\partial^*E$ and $\nu_E$ denote the reduced boundary
and the outer measure-theoretic unit normal of
$E\in\Cacc$, respectively.

\item $[E]$ and $\partial E$ denote the integral $3$-current
associated with $E\in\Cacc$ and its boundary integral
$2$-current, respectively, both with $\mathbb Z_2$-coefficients.

\item $\mathcal Z_n(M;\mathbb Z_2)$ denotes the space of
mod-$2$ $n$-cycles endowed with the flat topology, and
$\calF$ denotes the flat norm.

\item $\mathcal V_n(U)$ and $\mathcal{IV}_n(U)$ denote the
spaces of $n$-varifolds and integral $n$-varifolds in an
open set $U\subset M$, respectively.

\item $|K|$ denotes the multiplicity-one rectifiable
varifold associated with a two-rectifiable set $K$
of locally finite $\mathcal H_g^2$-measure.

\item $|\partial E|:=|\partial^*E|$ denotes the boundary
varifold associated with $E\in\Cacc$.

\item $\|V\|$ denotes the weight measure of a varifold $V$.

\item $\Theta^2(\|V\|,x)$, or simply $\Theta(V,x)$, denotes
the two-dimensional density of $\|V\|$ at $x$.

\item $T_x^r$ and $\eta_{x,r}$ denote the rescaled exponential
map and the blow-up map, respectively, defined by
\[
T_x^r:B_1(0)\subset T_xM\longrightarrow B_r(x),
\qquad T_x^r(z):=\exp_x(rz),
\]
and
\[
\eta_{x,r}:=r^{-1}\exp_x^{-1}
\quad\text{on }B_{\operatorname{inj}(M)}(x),
\]
for $0<r<\operatorname{inj}(M)$, so that
$\eta_{x,r}|_{B_r(x)}=(T_x^r)^{-1}$.

\item $\operatorname{Tan}(V,x)$ denotes the set of all
subsequential weak varifold limits of
$(\eta_{x,r})_\#V$ as $r\downarrow0$.

\item $\bfF_V$ denotes the varifold $\mathbf F$-metric
on the space of varifold, inducing the weak varifold topology
on subsets of bounded mass.

\item $\bfF$ denotes the $\mathbf F$-metric on $\Cacc$,
defined by
\begin{equation}
\bfF(E,F):=\Vol(E\triangle F)
+\bfF_V(|\partial E|,|\partial F|).
\label{eq:pairmetric}
\end{equation}
\end{itemize}

\subsection{Anisotropic energies and derivative conventions}
\label{subsec:anisotropic-derivatives}
We first state all assumptions on the anisotropic integrand and then fix the derivative conventions used in its first variation.

\vspace{5pt}
\noindent Throughout the paper, $\Phi\colon UM\to(0,\infty)$ is smooth and even in the second variable.If $(x,T)\in G_2(M)$, choose either unit normal $\nu_T$ to $T$ and set
\[
\Phi(x,T):=\Phi(x,\nu_T);
\]
this is independent of the choice of $\nu_T$ because $\Phi$ is even. When no confusion arises, we will often switch between these points of view without explicitly acknowledging it.  

By positive homogeneous extension of degree one, we define $\tilde{\Phi}$ on $TM\setminus\{0\}$ via
\begin{equation}
\tilde{\Phi}(x,v):=|v|_g\Phi\!\left(x,\frac{v}{|v|_g}\right),
\qquad (x,v)\in TM\setminus\{0\}.
\end{equation}

\noindent For $w,z\in T_xM$, its vertical Hessian at $(x,v)\in TM\setminus\{0\}$ is the symmetric bilinear form
\[
D_{vv}^2\tPhi(x,v)[w,z]
:=\left.\frac{\partial^2}{\partial s\,\partial t}\right|_{s=t=0}
\tPhi(x,v+sw+tz).
\]
We call $\Phi$ an \emph{elliptic integrand} if
\[
D_{vv}^2\tPhi(x,\nu)[w,w]>0
\]
 for every $(x,\nu)\in UM$ and all $w\neq 0\in\nu^\perp$.
By compactness of $UM$, any smooth elliptic integrand $\Phi$ satisfies the following quantitative bounds: there exists a constant $\lambda\geq 1$ such that
\begin{equation}\tag{$H_\lambda$}\label{eq:elliptic-class}
\lambda^{-1}\leq\Phi\leq\lambda,\quad
\|\Phi\|_{C^3(UM)}\leq\lambda,\quad
\lambda^{-1}|w|_g^2\leq D_{vv}^2\tilde{\Phi}(x,\nu)[w,w]\leq\lambda|w|_g^2
\end{equation}
for all $(x,\nu)\in UM$ and $w\in\nu^\perp$.

\vspace{5pt}
\noindent For a $2$-varifold $V$ and a Borel set $U\subset M$, we define its \textit{$\Phi$-anisotropic energy on $U$} by
\[
\Fphi(V;U):=\int_{G_2(U)}\Phi(x,T)\,\mathrm{d}V(x,T).
\]
We adopt the convention that omitting the second argument $U$ means taking the energy over the whole manifold $M$, i.e., $\Fphi(V):=\Fphi(V;M)$.
For a Caccioppoli set $\Omega\subset M$, we write
$\Fphi(\Omega;U)$ as shorthand for $\Fphi(|\partial\Omega|;U)$,
and for $T\in\mathcal{Z}_2(M;\mathbb{Z}_2)$, we write
$\Fphi(T;U)$ as shorthand for $\Fphi(|T|;U)$.

\vspace{10pt}
\noindent Having defined the anisotropic energy, we  now record its first variation .  If $F_t$ is the flow of a smooth vector field $X$, the \textit{anisotropic first variation along X} is defined by
\[
\delta_\Phi V(X)=\left.\frac{\dd}{\dd t}\right|_{t=0}\Fphi((F_t)_\#V).
\]

\noindent   For $H\geq0$,a varifold $V$ is said to have \emph{$H$-bounded
$\Phi$-first variation} if
\begin{equation}
 |\delta_\Phi V(X)|\leq H\int_U |X|_g\,\dd\|V\|,
 \qquad X\in C_c^1(TM).
 \label{eq:bounded-first-variation}
\end{equation}

\subsection{Anisotropic PMC hypersurfaces}
\label{subsec:anisotropic-pmc}

For $h\in C^\infty(M)$ and a Caccioppoli set $\Omega\in\Cacc$, define
\begin{equation}
\Aph(\Omega)=\Fphi(\partial\Omega)-\int_\Omega h\,\dd V_g.
\label{eq:Aph}
\end{equation}
This is the prescribed-mean-curvature analogue of the functional in
\cite[Section~2.3]{DDL24}.  If $\Sigma=\partial\Omega$ is smooth,
$\nu_\Omega$ is its outer unit normal, and $F_t$ is the flow of a
smooth ambient vector field $X$, its first variation is
\begin{equation}
\left.\frac{\dd}{\dd t}\right|_{t=0}\Aph(F_t(\Omega))
=\int_{\partial\Omega}
\bigl(H_\Phi[\partial\Omega,\nu_\Omega]-h\bigr)
\pair{X}{\nu_\Omega}\,\dd\mathcal H_g^2.
\label{eq:first-Aph}
\end{equation}
Here $H_\Phi[\partial\Omega,\nu_\Omega]$ denotes the anisotropic mean
curvature with respect to the outer unit normal.Consequently,if $\Omega$ is a
critical point of $\Aph$ , then
\[
H_\Phi[\partial\Omega,\nu_\Omega]=h|_{\partial\Omega}.
\]
we call
$\Sigma$ an \emph{$h$-anisotropic prescribed-mean-curvature boundary}
(or an \emph{$h$-APMC boundary}).More generally, we call a smooth immersion $\psi:S\to M$ with unit normal $\nu$ satisfying $H_\Phi[\psi,\nu]=h\circ\psi$
an \emph{$h$-anisotropic prescribed-mean-curvature surface}
(or an \emph{$h$-APMC surface}).

\vspace{8pt}
\noindent To give an explicit expression for the anisotropic mean curvature of a
smooth surface, we first fix the required covariant derivative notation.
 
\vspace{8pt}
\noindent We now introduce the covariant derivative notation.  Let $\nabla^g$ be the Levi--Civita connection of $g$; when there is no ambiguity, we write simply $\nabla=\nabla^g$.  If $\gamma:I\to M$ is a smooth curve and $W(t)\in T_{\gamma(t)}M$ is a smooth vector field along $\gamma$, its covariant derivative along $\gamma$ is denoted by
\[
\frac{D^gW}{\dd t}(t):=\nabla^g_{\dot\gamma(t)}W(t)
\in T_{\gamma(t)}M.
\]

Fix $(x,v)\in TM\setminus\{0\}$ and $Y,w\in T_xM$.  Choose a smooth curve $\gamma$ with $\gamma(0)=x$ and $\dot\gamma(0)=Y$, and let $V$ be the $g$-parallel vector field along $\gamma$ with $V(0)=v$.  The horizontal and vertical differentials of $\tPhi$ at $(x,v)$ are the covectors in $T_x^*M$ defined by
\[
d_x^{\mathrm h}\tPhi|_{(x,v)}(Y)
:=\left.\frac{\dd}{\dd t}\right|_{t=0}\tPhi(\gamma(t),V(t)),
\qquad
d_v\tPhi|_{(x,v)}(w)
:=\left.\frac{\dd}{\dd t}\right|_{t=0}\tPhi(x,v+tw).
\]
The first value depends only on $Y$, not on the chosen curve $\gamma$.  Their $g_x$-dual vectors
\[
D_x\tPhi(x,v),D_v\tPhi(x,v)\in T_xM
\]
are characterized by
\[
g_x(D_x\tPhi(x,v),Y)=d_x^{\mathrm h}\tPhi|_{(x,v)}(Y),
\qquad
g_x(D_v\tPhi(x,v),w)=d_v\tPhi|_{(x,v)}(w)
\]
for all $Y,w\in T_xM$.  Thus $D_x\tPhi$ always denotes the metric dual of the horizontal derivative, rather than a coordinate derivative. 

Using the same curve $\gamma$ and parallel vector field $V$, we define
the horizontal covariant derivative of $D_v\tPhi$ by
\[
\nabla_Y^{\mathrm h}[D_v\tPhi](x,v)
:=\left.\frac{D^g}{\dd t}\right|_{t=0}
\bigl(D_v\tPhi(\gamma(t),V(t))\bigr)\in T_xM.
\]
The differentiated expression is a vector field along $\gamma$, so
the derivative on the right is a covariant derivative along that
curve.  The fiber vector satisfies $D^gV/\dd t=0$; the value depends
only on $(x,v)$ and $Y$.
For a smooth two-sided surface $\Sigma$, we have
\begin{equation}
\delta_\Phi|\Sigma|(f\nu)=\int_\Sigma H_\Phi[\Sigma,\nu]f\,\dd\mathcal H_g^2,
\label{eq:H-convention}
\end{equation}
where $\nu$ and $H_\Phi[\Sigma,\nu](x)$ denote, respectively, the unit normal vector field on $\Sigma$ and \textit{anisotropic mean curvature} of $\Sigma$ at $x$ with respect to $\nu$.Under our conventions the anisotropic mean curvature admits the explicit form
\begin{equation}
H_\Phi[\Sigma,\nu]
=\operatorname{div}_\Sigma(D_v\tPhi(x,\nu))+
g_x(D_x\tPhi(x,\nu),\nu).
\label{eq:Hformula}
\end{equation}
Details of the derivation for this expression are given in Appendix~\ref{app:first-variation-formula}. We also note that the first‑variation formula for the anisotropic energy in the Euclidean setting may be found in \cite[Appendix~A]{DDG18}.

\vspace{8pt}

Let $F:S\to M$ be a smooth two-sided immersed $h$-APMC
surface with a global unit normal $\nu$ satisfying
$H_\Phi[F,\nu]=h\circ F$,we denote by $Q_{\Phi,h}^S$ the second variation of the
prescribed anisotropic functional, regarded as a symmetric
bilinear form on $C_c^\infty(S)$, where each
$u\in C_c^\infty(S)$ is identified with the normal
vector field $u\nu$.

\begin{definition}[Stability]
\label{def:strong-apmc-stability}
Let $U\subset M$ be open,a smooth two-sided $h$-APMC immersion $F:S\to M$, with chosen unit
normal $\nu$ satisfying $H_\Phi[F,\nu]=h\circ F$, is called
\emph{$h$-stable in $U$} if
\[
Q_{\Phi,h}^S(u\nu,u\nu)\ge0
\qquad\text{for every }u\in C_c^\infty(F^{-1}(U)).
\]

\end{definition}

For a smooth two-sided $h$-stable APMC surface $\Sigma$ in $U$,
the second-variation argument of \cite[Lemma~A.5, (A.13)]{DM17}gives the following estimate:
\begin{equation}
\int_\Sigma |A_\Sigma|^2\zeta^2\,\dd\mathcal H_g^2
\leq C\int_\Sigma
\bigl(|\nabla_\Sigma\zeta|^2+
(1+\|\nabla^g h\|_{C^0(U)})\zeta^2\bigr)
\,\dd\mathcal H_g^2,
\label{eq:pmc-stability-prelim}
\end{equation}
for every $\zeta\in C_c^1(\Sigma\cap U)$, where
$C=C(M,g,\lambda)$ and
$A_\Sigma(Y,Z)=g(\nabla_Y^g\nu,Z)$ is the second fundamental form.
For an immersion, the test functions and integrals are understood
on its parameter surface.

The h-APMC surface obtained via min-max theory may exhibit tangential self‑contact, and so we recall the following almost‑embedded immersion as introduced in \cite{ZZ19}.

\begin{definition}[Almost embedding]
Let $\phi\colon \Sigma \to M$ be a smooth immersion. We say that $\phi$ is an \emph{almost‑embedding} if, for any $p\in \phi(\Sigma)$ at which $\phi$ fails to be an embedding, there exists a neighborhood $U$ of $p$ such that
\begin{itemize}
    \item[(1)] $\Sigma \cap \phi^{-1}(U)$ decomposes into a disjoint union of finitely many pieces $\Sigma_i$, $i=1,\dots,l$, and $\phi$ restricted to each $\Sigma_i$ is an embedding;
    \item[(2)] for each $i$, every other image $\phi(\Sigma_j)$ with $j\neq i$ lies on one side of $\phi(\Sigma_i)$.
\end{itemize}
\end{definition}
\noindent We identify $\Sigma$ with its image and refer to the image as \emph{almost‑embedded}. Moreover, we let $\mathcal{S}(\Sigma)$ denote the set of points in $\Sigma$ where the immersion fails to be an embedding, and let $\mathcal{R}(\Sigma)$ denote $\Sigma \setminus \mathcal{S}(\Sigma)$.

\section{Generic prescription functions and touching sets}
\label{sec:good-local}

In this section,we first prove that prescription functions satisfying the required nondegeneracy along their zero sets form an open and dense class; we then analyze touches between APMC surface whose anisotropic mean curvature is given by such functions.

\subsection{\texorpdfstring{$\Phi$}{Phi}-good prescription functions: open and dense}
\begin{definition}
For a smooth $h$ having $0$ as a regular value, write
$Z(h)=h^{-1}(0)$, $\nu_h=\nabla h/|\nabla h|$, and
$K_{\Phi,h}=H_\Phi[Z(h),\nu_h]$.
A smooth $h$ is $\Phi$-good if it is Morse, $0$ is a regular value, and at every $x\in Z(h)$ with $K_{\Phi,h}(x)=0$ one has
$d(K_{\Phi,h})_x\in T_x^*Z(h)\setminus\{0\}$.
\end{definition}

The role of this condition is to exclude infinite-order touching between
APMC sheets with opposite outer normals .  Such a contact would force
the anisotropic mean curvature of the zero level set, together with all of
its tangential derivatives, to vanish at the contact point, contrary
to the defining condition on $K_{\Phi,h}$.  Once infinite-order contact
is excluded, a vanishing-order stratification gives the touching-set
estimate for almost embedded APMC boundaries proved in
Theorem~\ref{thm:touching} and
Corollary~\ref{cor:apmc-boundary-touching}.  This estimate is a local
ingredient in the regularity theory for the APMC
min--max construction.

This condition plays the same role as the special class of
prescription functions used by Zhou--Zhu in their isotropic PMC min--max
theory: they consider Morse functions with regular zero level set whose
mean curvature vanishes only to finite order
\cite[Proposition~3.7]{ZZ20}.  For $\Phi\equiv1$, our condition
imposes the stronger requirement that every zero of that mean
curvature be nondegenerate as a function on the zero level.
Zhou--Zhu establish density of their class by evolving the zero
level by mean curvature flow, using short-time existence and unique
continuation for the curvature evolution equation
\cite[proof of Proposition~3.7 and Section~3.3]{ZZ20}.

Here we prove $C^\infty$-density of $\Phi$-good prescription functions by a
transversality argument.  We express the anisotropic mean curvature
of a regular level set as a function of the two-jet of its defining
function, prove the relevant jet-bundle map is a submersion, and
apply jet transversality to obtain the curvature condition.

\vspace{8pt}
We begin by recalling the definition of jet-bundle.
For an integer $k\geq0$ and $x\in M$, let $f$ and $g$ be smooth
real-valued functions defined on neighborhoods of $x$.  We say that $f$
and $g$ have the same $k$-jet at $x$ if, for one (and hence every)
coordinate chart $\psi$ about $x$,
\[
\partial^\alpha(f\circ\psi^{-1})(\psi(x))
=\partial^\alpha(g\circ\psi^{-1})(\psi(x))
\qquad\text{for every multi-index $\alpha$ with $|\alpha|\leq k$}.
\]
The resulting equivalence class $j_x^kf$ is the $k$-jet of $f$ at $x$.
The fiber $J_x^k(M,\R)$ is the set of all $k$-jets at $x$, and
\[
J^k(M,\R)=\bigsqcup_{x\in M}J_x^k(M,\R),
\qquad
\pi_k(j_x^kf)=x,
\]
is the $k$-jet bundle.  Every $f\in C^\infty(M)$ determines a smooth
$k$-jet by
\[
j^kf:M\longrightarrow J^k(M,\R),
\qquad x\longmapsto j_x^kf.
\]
For $k=2$, the Levi--Civita connection of $g$ gives the fiberwise
identification
\[
J_x^2(M,\R)
\cong\R\oplus T_x^*M\oplus\operatorname{Sym}^2T_x^*M,
\qquad
j_x^2f\longmapsto\bigl(f(x),df_x,(\nabla^2f)_x\bigr).
\]
Thus a two-jet will be written as $z=(x,r,p,Q)$, with
$p\in T_x^*M$ and $Q\in\operatorname{Sym}^2T_x^*M$.

\vspace{8pt}
We next express the anisotropic mean curvature of a regular level set in terms of this two-jet and then extend the expression to the open set of all two-jets with nonzero first-order part.  We use the horizontal and vertical derivative notation fixed in \cref{subsec:anisotropic-pmc}.

\vspace{5pt}
We first collect the terms in the expression of anisotropic mean curvature of a regular level set depending only on the base point and the chosen unit normal.  For $(x,\nu)\in UM$ and an orthonormal basis $e_1,e_2$ of $\nu^\perp$, set
\begin{equation}
b_\Phi(x,\nu)=
\sum_{\alpha=1}^2
g_x\bigl(\nabla_{e_\alpha}^{\mathrm h}[D_v\tPhi](x,\nu),e_\alpha\bigr)
+d_x^{\mathrm h}\tPhi|_{(x,\nu)}(\nu).
\label{eq:b-Phi}
\end{equation}
The first term is a trace on $\nu^\perp$, so this definition is independent of the chosen orthonormal basis and gives a smooth function $b_\Phi:UM\to\R$.

Let $f\in C^\infty(M)$ and let $x\in M$ satisfy $df_x\neq0$.  Near $x$, the level set
\[
\Sigma_{f,x}=f^{-1}(f(x))
\]
is a smooth hypersurface.  Put
\[
\nu_f(x)=\frac{(df_x)^{\sharp_g}}{|df_x|_g}.
\]
Let
\[
P_f(x):T_xM\longrightarrow T_x\Sigma_{f,x}
=\ker df_x=\nu_f(x)^\perp
\]
be the $g_x$-orthogonal projection.
If $e_1,e_2$ is an orthonormal basis of $T_x\Sigma_{f,x}$, differentiation of $\nu_f$ gives
\begin{equation}
\nabla_{e_\alpha}\nu_f
=\frac1{|df_x|_g}P_f(x)\bigl((\nabla^2f)_x^{\sharp_g}e_\alpha\bigr),
\qquad \alpha=1,2.
\label{eq:level-normal-derivative}
\end{equation}
The chain rule in the horizontal--vertical splitting gives
\[
\nabla_{e_\alpha}\bigl(D_v\tPhi(\,\cdot\,,\nu_f(\,\cdot\,))\bigr)
=\nabla_{e_\alpha}^{\mathrm h}[D_v\tPhi](x,\nu_f(x))
+A_\Phi(x,\nu_f(x))\nabla_{e_\alpha}\nu_f.
\]
Substituting this identity and \eqref{eq:level-normal-derivative} into \eqref{eq:Hformula} yields
\begin{align}
H_\Phi[\Sigma_{f,x},\nu_f](x)
&=b_\Phi(x,\nu_f(x)) \notag\\*
&\quad+
\frac1{|df_x|_g}
\tr_{T_x\Sigma_{f,x}}\!\left(
A_\Phi(x,\nu_f(x))\circ P_f(x)\circ
(\nabla^2f)_x^{\sharp_g}\big|_{T_x\Sigma_{f,x}}
\right).
\label{eq:classical-level-H}
\end{align}

We now extend this expression to the two-jet bundle.  Let
\[
J^2_{\rm reg}(M,\R)=\{(x,r,p,Q):p\neq0\}.
\]
For $p\neq0$, define
\[
\nu(p)=\frac{p^{\sharp_g}}{|p|_g},\qquad
T(p)=\ker p,\qquad
P(p)=\operatorname{Id}_{T_xM}-\nu(p)\otimes\nu(p)^{\flat_g}.
\]
For $Q\in\operatorname{Sym}^2T_x^*M$, define
\begin{equation}
\mathcal H_\Phi(x,r,p,Q)
=b_\Phi(x,\nu(p))+
\frac1{|p|_g}\tr_{T(p)}\!\left(
A_\Phi(x,\nu(p))\circ P(p)\circ Q^{\sharp_g}|_{T(p)}
\right).
\label{eq:two-jet-H}
\end{equation}
Formula \eqref{eq:classical-level-H} states precisely that
\[
\mathcal H_\Phi(j_x^2f)
=H_\Phi[f^{-1}(f(x)),\nu_f](x)
\qquad\text{whenever }df_x\neq0.
\]
 To keep track simultaneously of the level value and its anisotropic mean curvature, define
\[
\Theta_\Phi:J^2_{\rm reg}(M,\R)\longrightarrow\R^2,
\qquad
\Theta_\Phi(x,r,p,Q)=\bigl(r,\mathcal H_\Phi(x,r,p,Q)\bigr).
\]

The following submersion property provides the key differential-topological ingredient for the jet-transversality argument.
\begin{lemma}[Jet-bundle submersion]
\label{lem:two-jet}
The function $\Theta_\Phi$ is smooth and is a submersion at every point of $J^2_{\rm reg}(M,\R)$.
\end{lemma}

\begin{proof}
The maps $p\mapsto |p|_g^{-1}$, $p\mapsto\nu(p)$, and $p\mapsto P(p)$ are smooth on $T^*M\setminus\{0\}$.  The trace in \eqref{eq:two-jet-H} is the contraction of smooth tensor fields and is linear in $Q$.  Since $b_\Phi$ and $A_\Phi$ are smooth, both $\mathcal H_\Phi$ and $\Theta_\Phi$ are smooth.

Fix $z=(x,r,p,Q)\in J^2_{\rm reg}(M,\R)$.  For
$B\in\operatorname{Sym}^2T_x^*M$, vary only the Hessian component along the fiberwise path
\[
z_B(t)=(x,r,p,Q+tB).
\]
Differentiating \eqref{eq:two-jet-H} at $t=0$ gives
\begin{equation}
d\mathcal H_\Phi|_z\bigl(0,0,0,B\bigr)
=\frac1{|p|_g}\tr_{T(p)}\!\left(
A_\Phi(x,\nu(p))\circ P(p)\circ B^{\sharp_g}|_{T(p)}
\right).
\label{eq:H-A-derivative}
\end{equation}
Take $B=g_x$.  Then $B^{\sharp_g}=\operatorname{Id}_{T_xM}$ and $P(p)|_{T(p)}=\operatorname{Id}_{T(p)}$.  If $e_1,e_2$ is an orthonormal basis of $T(p)$, the ellipticity assumption \eqref{eq:elliptic-class} implies
\[
d\mathcal H_\Phi|_z(0,0,0,g_x)
=\frac1{|p|_g}\sum_{\alpha=1}^2
D_{vv}^2\tPhi(x,\nu(p))[e_\alpha,e_\alpha]
\geq\frac{2}{\lambda|p|_g}>0.
\]
On the other hand, varying only the zero-order component gives
\[
d\Theta_\Phi|_z(0,1,0,0)=(1,0),
\]
because \eqref{eq:two-jet-H} is independent of $r$.Therefore $d\Theta_\Phi|_z$ is surjective.  Since $z$ was arbitrary, $\Theta_\Phi$ is a submersion everywhere.
\end{proof}

The submersion established in Lemma~3.2 allows us to apply jet transversality to prove the density of $\Phi$-good functions; we also establish their openness in the following proposition.
\begin{proposition}
\label{thm:good-generic}
The $\Phi$-good functions form a $C^3$-open and $C^\infty$-dense subset of $C^\infty(M)$. 
\end{proposition}

\begin{proof}
Since $M$ is compact, every Morse function has only finitely many critical points. Therefore, we may add an arbitrarily small constant $c$ so that $0$ becomes a regular value of $f+c$. Combined with the density and $C^2$-openness of Morse functions, this yields

\[
\mathscr M_0
=\{h\in C^\infty(M):h\text{ is Morse and }0\text{ is a regular value of }h\}
\]
is $C^2$-open and $C^\infty$-dense.

Openness follows from (3.3): the curvature $K_{\Phi,h}$ and its
tangential differential along $Z(h)$ depend continuously on the
derivatives of $h$ of order at most three. Since $Z(h)$ is compact,
the condition
\[
|K_{\Phi,h}|+|d_{Z(h)}K_{\Phi,h}|>0
\qquad\text{on }Z(h)
\]
has a uniform positive lower bound and therefore persists under
sufficiently small $C^3$ perturbations. The Morse condition and
the regularity of the zero level are also preserved under such
perturbations. Hence the set of $\Phi$-good functions is
$C^3$-open.

It remains only to prove density of the curvature condition.  By
Lemma~\ref{lem:two-jet},
\begin{equation}
\mathscr W_\Phi
=\Theta_\Phi^{-1}(0,0)
=\{r=0,\ \mathcal H_\Phi=0\}
\subset J^2_{\rm reg}(M,\R)
\label{eq:bad-jet}
\end{equation}
is a smooth codimension-two submanifold.  Since
$J^2_{\rm reg}(M,\R)$ is open in $J^2(M,\R)$, $\mathscr W_\Phi$ is
also an embedded, not necessarily closed, submanifold of
$J^2(M,\R)$.

Suppose that $0$ is a regular value of $h$.  On $Z(h)=h^{-1}(0)$ define
\[
\nu_h=\frac{\nabla h}{|\nabla h|},
\qquad
K_{\Phi,h}=H_\Phi[Z(h),\nu_h],
\]
and on the open set $\{dh\neq0\}$ define
\[
\widehat K_{\Phi,h}(x)=\mathcal H_\Phi(j_x^2h).
\]
The restriction of $\widehat K_{\Phi,h}$ to $Z(h)$ is $K_{\Phi,h}$,and we have

\[
\Theta_\Phi\circ j^2h=(h,\widehat K_{\Phi,h})
\]
Suppose $x\in Z(h)$ and $K_{\Phi,h}(x)=0$, Since $\mathscr W_\Phi=\Theta_\Phi^{-1}(0,0)$ and
$d\Theta_\Phi$ is surjective, $j^2h$ is transverse to
$\mathscr W_\Phi$ at $x$ if and only if
\begin{equation}
\operatorname{rank}d(h,\widehat K_{\Phi,h})_x=2.
\label{eq:rank-composite}
\end{equation}
Indeed, at $z=j_x^2h$ one has
$T_z\mathscr W_\Phi=\ker d\Theta_\Phi|_z$.  Thus transversality
implies the surjectivity of $d(\Theta_\Phi\circ j^2h)|_x$.
Conversely, if the latter differential is surjective, then for every
$v\in T_zJ^2_{\rm reg}(M,\R)$ one can choose $X\in T_xM$ such that
$d\Theta_\Phi|_z(v-d(j^2h)|_xX)=0$; hence
$v\in d(j^2h)|_x(T_xM)+T_z\mathscr W_\Phi$, which is transversality.

Because $dh_x\neq0$ and $T_xZ(h)=\ker dh_x$, elementary linear algebra
gives
\begin{equation}
\operatorname{rank}d(h,\widehat K_{\Phi,h})_x=2
\quad\Longleftrightarrow\quad
\left.d\widehat K_{\Phi,h}|_x\right|_{\ker dh_x}\neq0
\quad\Longleftrightarrow\quad
d(K_{\Phi,h})_x\neq0\text{ in }T_x^*Z(h).
\label{eq:rank-two}
\end{equation}
Consequently, for $h\in\mathscr M_0$, the curvature condition is
equivalent to
\begin{equation}
j^2h\pitchfork\mathscr W_\Phi.
\label{eq:good-transversality}
\end{equation}
By the jet transversality theorem
\cite[Chapter~3, \S2, Theorem~2.8]{Hir76}, the functions
$h\in C^\infty(M)$ such that $j^2h\pitchfork\mathscr W_\Phi$ form a
residual set, and hence dense subset of $C^\infty(M)$.  Since
$\mathscr M_0$ is open and dense, its intersection with this dense set
is dense, and every function in the intersection is $\Phi$-good.

\end{proof}

\subsection{Touching set estimates}
We establish a touching-set estimate for APMC surfaces with a
$\Phi$-good prescription and record its unique continuation.  The following theorem
is the anisotropic counterpart of \cite[Theorem~3.10]{ZZ20}.

\begin{theorem}[Touching structure]
\label{thm:touching}
Let $h$ be $\Phi$-good and let $U\subset M$ be a connected open set.
Let $\Sigma_1,\Sigma_2\subset U$ be connected, smooth, properly
embedded surfaces without boundary in $U$, equipped with smooth unit
normals $\nu_i$ satisfying
\[
H_\Phi[\Sigma_i,\nu_i]=h|_{\Sigma_i},\qquad i=1,2.
\]
Assume that they lie locally to one side of one another at every
intersection point.  Then either $\Sigma_1=\Sigma_2$, or their touching set
$\Sigma_1\cap\Sigma_2$ is contained in a countable union of connected,
smoothly embedded one-dimensional submanifolds of $U$.
Moreover, in the latter case, $T$ has locally finite
one-dimensional Hausdorff measure. More precisely,
\[
\mathcal H_g^1(T\cap K)<\infty
\qquad\text{for every compact }K \Subset U.
\]

\end{theorem}

\begin{proof}
Fix a touching point $p\in\Sigma_1\cap\Sigma_2$.  Let $(y,z)$ denote geodesic normal coordinates centered at $p$, $y=(y_1,y_2)$, for which $\partial_{y_1}|_0,\partial_{y_2}|_0$ span the common tangent plane $T_p\Sigma_1=T_p\Sigma_2$ and $\nu_1(p)=\partial_z|_0$.
In a sufficiently small coordinate
cylinder, the two surfaces are graphs $z=u_i(y)$ over a common
connected disk $D\subset\R^2$, with
\[
u_1(0)=u_2(0)=0,\qquad Du_1(0)=Du_2(0)=0.
\]
Without loss of generality, we assume that $u_1\geq u_2$ throughout $D$.

Let $H_\Phi[u]$ denote the anisotropic mean curvature of
the graph $z=u(y)$ with respect to its upward unit normal, namely
the unit normal having positive $g$-inner product with $\partial_z$.
According to the direction of $\nu_2(p)$, the graph equations are therefore
\[
\begin{aligned}
H_\Phi[u_1]&=h(y,u_1(y)),\\
H_\Phi[u_2]&=
\begin{cases}
h(y,u_2(y)),&\text{if }\nu_2(p)=\nu_1(p),\\
-h(y,u_2(y)),&\text{if }\nu_2(p)=-\nu_1(p).
\end{cases}
\end{aligned}
\]
For $w=u_1-u_2$,we have
\begin{equation}
L_{u_2,u_1}w
=H_\Phi[u_1]-H_\Phi[u_2].
\label{eq:secant}
\end{equation}
Here, $L$ is a linear second-order operator with
smooth coefficients and a uniformly negative definite principal
part.

\vspace{8pt}
\noindent If $\nu_2(p)=\nu_1(p)$, the difference satisfies the homogeneous
equation
\begin{equation}
\bigl(L_{u_2,u_1}-c_h\bigr)w=0,
\qquad
c_h(y):=\int_0^1
\partial_z h\bigl(y,u_2(y)+s w(y)\bigr)\dd s.
\end{equation}
If $\nu_2(p)=-\nu_1(p)$, it instead satisfies
\begin{equation}
L_{u_2,u_1}w=h(y,u_1(y))+h(y,u_2(y)).
\label{eq:opposite-equation}
\end{equation}
The following same-orientation comparison follows from Serrin's
strong maximum principle \cite{Ser70}.

\begin{lemma}[Maximum principle for APMC sheets]
\label{lem:apmc-same-orientation}
Let $h\in C^\infty(U)$.  Let $\Sigma_1,\Sigma_2\subset U$ be
connected, smooth, properly embedded h-APMC surfaces without boundary
in $U$, with chosen unit normals $\nu_i$, and suppose that they lie
locally to one side of one another at each intersection point.
If $p\in\Sigma_1\cap\Sigma_2$ and $\nu_1(p)=\nu_2(p)$, then
$\Sigma_1=\Sigma_2$, with the same chosen normals.
\end{lemma}

Returning to the theorem, assume that $\Sigma_1\neq\Sigma_2$.
By Lemma~\ref{lem:apmc-same-orientation}, the chosen normals are
opposite at every touching point.  It remains to prove a local
curve covering near each such point; a countable covering by
coordinate cylinders will then give the stated conclusion.
Let
\[
T=\{y\in D:w(y)=0\}.
\]
By our assumption on $\Sigma_1$ and $\Sigma_2$,every point of $T$ also satisfies $Dw=0$.

\emph{Finite-order touch.}
We use the vanishing-order stratification and implicit-function
argument from the proof of \cite[Lemma~1.9]{HS89}.
For each integer $k\geq2$, let
\[
T_k=\bigl\{y\in T:D^\alpha w(y)=0\text{ for }|\alpha|<k,
\quad D^\beta w(y)\neq0\text{ for some }|\beta|=k\bigr\}.
\]
For $y_0\in T_k$, choose a multi-index $\alpha$ with $|\alpha|=k-1$
and an index $j\in\{1,2\}$ such that
$\partial_jD^\alpha w(y_0)\neq0$.  The implicit function theorem
gives a neighborhood $V$ of $y_0$ on which
\[
\Gamma=\{y\in V:D^\alpha w(y)=0\}
\]
is a connected smooth curve.  Thus we have
$T_k\cap V\subset\Gamma$.Second countability gives a countable
subcover of each $T_k$ by such neighborhoods.  Taking the union over
$k$ covers all finite-order vanishing points by countably many smooth
curves.

\emph{Exclusion of infinite-order touch.}
Suppose all derivatives of $w$ vanish at $0$. Write
$h_i(y)=h(y,u_i(y))$ for $i=1,2$. The identity
\[
h_1(y)-h_2(y)
=w(y)\int_0^1\partial_z h\bigl(y,u_2(y)+s w(y)\bigr)\,\dd s
\]
shows that $h_1-h_2$ vanishes to infinite order at $0$.
By \eqref{eq:opposite-equation}, their sum is
$ L_{u_2,u_1}w$, which also vanishes to infinite order. Thus
each $h_i$ vanishes to infinite order at $0$,in particular we have $h(p)=0$.

Differentiating
$h_1$ at $0$ gives
\[
\partial_{y_i}h(p)+\partial_z h(p)D_i u_1(0)=0,
\qquad i=1,2.
\]
Since $Du_1(0)=0$ we have $\partial_{y_i}h(p)=0$ and consequently  we have $\partial_z h(p)\neq0$, otherwise it would contradict with the fact  $0$ is a regular value of $h$.
The implicit function theorem therefore represents $Z(h)$ near $p$
as a graph $z=v_0(y)$.Then we have
\[
\begin{aligned}
h(y,z)&=a(y,z)\bigl(z-v_0(y)\bigr),\\
a(y,z)&=\int_0^1\partial_z h\bigl(y,v_0(y)+s(z-v_0(y))\bigr)\,\dd s.
\end{aligned}
\]
Here $a$ is smooth and $a(p)=\partial_z h(p)\neq0$. After shrinking the cylinder to a smaller neighbourhood of $p$, we may assume $a$ is nowhere vanishing on this neighbourhood.
and hence
\[
u_i(y)-v_0(y)=\frac{h_i(y)}{a(y,u_i(y))},
\qquad i=1,2,
\]
vanishes to infinite order at $0$.

Combined with the fact that $H_\Phi[u_1] = h_1$ vanishes to infinite order at $0$, and that the curvature $H_\Phi[v_0]$ together with its tangential differential along $Z(h)$ at $p$ depends only on the jet of order at most three of $h$ at $p$.We deduce that
\[
H_\Phi[v_0](p)=0,
\qquad d(H_\Phi[v_0])_p=0\quad\text{in }T_p^*Z(h),
\]
contrary to the definition of $\Phi$-good function.Thus every point of $T$ has finite vanishing order.

To see $T$ has locally finite
one-dimensional Hausdorff measure,fix $y_0\in T_k$ ,since $w$ has vanishing order $k$
at $y_0$,we may
choose two linearly independent unit vectors
$v_1,v_2$ such that
$\partial_{v_i}^k w(y_0)\neq 0$
for $i=1,2$ and these derivatives remain nonzero in
some disk $B_r(y_0) \subset D$. Rolle's theorem then
implies that every line parallel to $v_i$ meets
$T\cap B_r(y_0)$ in at most $k$ points.

Let $\pi_i:\mathbb R^2\to v_i^\perp$ denote the
orthogonal projection. We denote by $J^T\pi_i(y)$ the one-dimensional
tangential Jacobian of $\pi_i$ restricted to
$\operatorname{Tan}^1(T,y)$.
The area formula gives
\[
\begin{aligned}
\int_{T\cap B_r(y_0)} J^T\pi_i\,d\mathcal H^1
&=
\int_{v_i^\perp}
\mathcal H^0\bigl(
T\cap B_r(y_0)\cap\pi_i^{-1}(s)
\bigr)\,d\mathcal H^1(s)\\
&\le 2kr,\qquad i=1,2.
\end{aligned}
\]
Since $v_1$ and $v_2$ are linearly independent,
there exists $c>0$ such that
\[
J^T\pi_1+J^T\pi_2\ge c
\]
for $\mathcal H^1$-almost every point of $T$.
Combining the preceding estimates for $i=1,2$ with the lower bound for the tangential Jacobians,we obtain
\[
\mathcal H^1(T\cap B_r(y_0))
\le \frac{4kr}{c}<\infty.
\]

A finite covering argument in local graph charts
yields $\mathcal H_g^1(T\cap K)<\infty$
for every compact $K\Subset U$.

\end{proof}
We record the following immediate consequences of the touching
analysis in Theorem~\ref{thm:touching}.

\begin{corollary}[Touching sets of almost embedded APMC surfaces]
\label{cor:apmc-boundary-touching}
Let $h$ be $\Phi$-good and let $U\subset M$ be open.  If
$\Sigma$ is a smooth almost embedded surface
in $U$,whose anisotropic mean curvature is $h$ with respect to its outer unit normal, then its touching set
$\mathcal S(\Sigma)$ is contained in a countable union of connected,
smoothly embedded one-dimensional submanifolds of $U$.
In particular, $\mathcal S(\Sigma)$ is countably $1$-rectifiable,
$\mathcal H_g^2(\mathcal S(\Sigma))=0$, and
$\mathcal R(\Sigma)=\Sigma\setminus\mathcal S(\Sigma)$ is open and
dense in $\Sigma$.
\end{corollary}

\begin{corollary}[Unique continuation for immersed APMC surfaces]
\label{cor:apmc-unique-continuation}
Let $h$ be $\Phi$-good, and let $\Sigma_1,\Sigma_2$ be two connected,
smooth, properly immersed surfaces without boundary, of prescribed
anisotropic mean curvature $h$, in a connected open subset
$U\subset M$.  If $\Sigma_1$ and $\Sigma_2$ coincide on a nonempty
open neighborhood $U'\subset U$ with $\Sigma_1\cap U'\neq\varnothing$,
then $\Sigma_1=\Sigma_2$ as subsets of $U$.
\end{corollary}
\section{Prescribed anisotropic min--max in dimension three}
\label{sec:transfer}

In this section,we establish the min-max theory of anisotropic prescribed-mean-curvature surface.Throughout this section $(M^3,g)$ is closed,$\Phi$ satisfies
\eqref{eq:elliptic-class} and $h$ is a $\Phi$-good function.The conclusion
is the prescribed-mean-curvature counterpart, in dimension three, of
\cite[Theorem~B]{DDL24}.

\subsection{ONVP sweepouts and the min--max construction}
We first give the following definition:

\begin{definition}[Sweepout]
A continuous map $\Omega:[0,1]\to\Cacc$ is called a \emph{sweepout of $M$} if
\begin{enumerate}[label=(\roman*)]
\item $\Omega(0)=\varnothing$ and $\Omega(1)=M$ ;
\item $\Omega$ has no concentration of mass,i.e
\begin{equation}
\lim_{r\downarrow0}\sup_{t\in[0,1]}\sup_{p\in M}
\Fphi(\Omega(t);B_r(p))=0.
\label{eq:no-concentration}
\end{equation}
\end{enumerate}
We denote the collection of all sweepouts by $\Pi(M)$,we then define the $\Wph$-width with respect to $\Phi$ and $h$.
\end{definition}

\begin{definition}[First prescribed width]
For $h\in C^\infty(M)$, define
\begin{equation}
\Wph:=\inf_{\Omega\in\Pi(M)}\sup_{t\in[0,1]}\Aph(\Omega(t)),
\qquad
\Wphi:=\mathbf W_{\Phi,0}.
\label{eq:width}
\end{equation}
\end{definition}

We first establilsh a lower bound for the prescribed width.

\begin{lemma}[Strict endpoint gap]\label{lem:strict-endpoint-gap}
For every $h\in C^\infty(M)$,
\begin{equation}
\max\left\{0,-\int_M h\,dV_g\right\}
<
\mathbf W_{\Phi,h}
<\infty.
\label{eq:strict-endpoint-gap}
\end{equation}
\end{lemma}

\begin{proof}
By \eqref{eq:elliptic-class} and the isoperimetric inequality
\cite[Theorem~2.3]{ZZ20}, there exist $c>0$ and $v_0>0$ such that
\[
\Fphi(E)\ge c\,\min\{\Vol(E),\Vol(M\setminus E)\}^{2/3}
\]
whenever $\min\{\Vol(E),\Vol(M\setminus E)\}$ is at most $v_0$.

Let $H=\|h\|_{C^0(M)}$ and choose
\[
0<v<\min\left\{v_0,\frac12\Vol(M)\right\}
\]
so that
\[
cv^{2/3}-Hv>0.
\]
Every sweepout passes through slices of volume $v$ and
$\Vol(M)-v$. Hence, for some $t_-,t_+\in(0,1)$,
\[
\mathcal A_{\Phi,h}(\Omega(t_-))
\ge cv^{2/3}-Hv>0,
\]
and
\[
\mathcal A_{\Phi,h}(\Omega(t_+))
\ge
-\int_M h\,dV_g+cv^{2/3}-Hv.
\]
Therefore
\[
\sup_{t\in[0,1]}\mathcal A_{\Phi,h}(\Omega(t))
>
\max\left\{0,-\int_M h\,dV_g\right\}.
\]
Taking the infimum over all sweepouts gives the desired lower bound.
Finally, a fixed smooth sweepout with uniformly bounded anisotropic
energy gives $\mathbf W_{\Phi,h}<\infty$.
\end{proof}
Following the ONVP sweepout framework introduced by
Chodosh--Liokumovich--Spolaor \cite{CLS22}, and subsequently
adapted to the anisotropic min--max setting in \cite{DDL24},
we introduce its prescribed anisotropic counterpart for the functional
$\mathcal A_{\Phi,h}$.
$\mathcal A_{\Phi,h}$.
\begin{definition}[Optimal nested volume-parameterized sweepout]
A sweepout $\Omega\in\Pi(M)$ is an \emph{optimal nested volume-parameterized} (ONVP) sweepout for $(\Phi,h)$ if
\begin{enumerate}[label=(\roman*)]
\item \emph{Optimal:}
\[
\sup_{t\in[0,1]}\Aph(\Omega(t))=\Wph;
\]
\item \emph{Nested:}
\[
\Omega(s)\subset\Omega(t)\quad\text{modulo null sets}
\qquad(0\le s\le t\le1);
\]
\item \emph{Volume-parameterized:}
\[
\Vol(\Omega(t))=t\Vol(M),\qquad t\in[0,1].
\]
\end{enumerate}
\end{definition}
We will also use the following notion of mass-ratio upper bound.
\begin{definition}[Mass-ratio bound]
For $V\in\mathcal V(M)$ and $r_0>0$, set
\begin{equation}
\operatorname{MR}_{\Phi,r_0}(V)
:=\sup_{p\in M}\sup_{0<r<r_0}
\frac{\mathbf{F}_{\Phi}(V;B_r(p))}{r^2}.
\label{eq:mass-ratio}
\end{equation}
For $E\in\Cacc$, set
\[
\operatorname{MR}_{\Phi,r_0}(E)
:=\operatorname{MR}_{\Phi,r_0}(|\partial E|).
\]
We say that $V$ (respectively $E$) has \emph{mass-ratio bound $C$ below scale $r_0$} if
$\operatorname{MR}_{\Phi,r_0}(V)\le C$ (respectively
$\operatorname{MR}_{\Phi,r_0}(E)\le C$).
\end{definition}
The mass-ratio control established in the next theorem is a key ingredient in the min--max existence theory for anisotropic prescribed mean curvature hypersurfaces.
\begin{theorem}[ONVP sweepouts with controlled mass ratio]
\label{thm:ONVP}
\label{thm:MR-ONVP}
Suppose $h$ is a smooth function,$\Phi$ be an elliptic integrand satisfying
\eqref{eq:elliptic-class}and $\|h\|_{C^0(M)}\le\Lambda$,then there is an ONVP sweepout with
 $\operatorname{MR}_{\Phi,r_0}(E)\leq C$ for any $r_0\in(0,\frac{inj(M)}{2})$,
where $C$ depend only on $M,g,\lambda$ and $\Lambda$.
\end{theorem}
\begin{proof}
The proof of \cite[Theorem 3.5]{DDL24} uses the volume term
$c\,\operatorname{Vol}(\Omega)$ only through its Lipschitz control with
respect to the symmetric-difference distance. Namely,
\[
\left|
c\,\operatorname{Vol}(E)-c\,\operatorname{Vol}(F)
\right|
\leq
c\,\operatorname{Vol}(E\triangle F).
\]
For the $h$-weighted volume term we have, analogously,
\[
\left|
\int_E h\,d\operatorname{Vol}
-
\int_F h\,d\operatorname{Vol}
\right|
\leq
\|h\|_{C^0}\operatorname{Vol}(E\triangle F).
\]
Therefore, the same proof yields the existence of an ONVP sweepout for
$\mathcal A_{\Phi,h}$.

We now follow the proof of \cite[Theorem 3.9]{DDL24}. We use the same triangulations and notation for their simplices as in
that proof .Throughout that
proof, we replace $\Phi^c$ by $\mathcal A_{\Phi,h}$, every term of the form$c\,\operatorname{Vol}(E)$
by$\int_E h\,d\operatorname{Vol}$,
and the upper bound $\Lambda$ for $c$ by an upper bound $\Lambda$ for
$\|h\|_{C^0}$.

Except for the quantitative estimates involving the volume term, all the other parts of the proof of \cite[Theorem 3.9]{DDL24} carry over without change.
It therefore remains only to verify the corresponding estimates for the $h$-weighted volume term.
Whenever a term of the form $c\,\Vol(E)$ appears in the estimates of
\cite[Theorem 3.9]{DDL24}, the corresponding term in the present setting
is estimated by $\left|\int_E h\,\dd V_g\right|
\leq \Lambda\Vol(E)$. For each simplex $\sigma$ at scale $r$, one has
$\Vol(\sigma)\leq Cr^3$, and hence the corresponding $h$-weighted
volume term is bounded by $C\Lambda r^3$.Since the anisotropic energy terms are of order $r^2$,
choosing the initial triangulation scale sufficiently small allows
these $r^3$-terms to be absorbed into the same $r^2$ bounds used in
\cite[Theorem 3.9]{DDL24}.Thus, after possibly decreasing the initial triangulation scale and
enlarging the constants in the anisotropic energy bounds of order
$r^2$, with the new constants depending additionally only on
$\Lambda$, all the quantitative inequalities in the proof remain
valid.This completes the verification of the quantitative estimates that differ
from those in the proof of \cite[Theorem 3.9]{DDL24}.

Therefore we obtain an
ONVP sweepout $\{\Omega(t)\}_{t\in[0,1]}$ and constants $C>0$ and
$r_0>0$, depending only on $M$, $g$, $\lambda$, and an upper bound $\Lambda$ for
$\|h\|_{C^0}$, such that
\[
\Fphi(\Omega(t);B_r(p))
\leq
Cr^2
\]
for every $t\in[0,1]$, $p\in M$, and $0<r<r_0$,from the hypothesis \eqref{eq:elliptic-class}.
\end{proof}
We next record the min--max sequences and their critical pairs associated with the width $\mathbf W_{\Phi,h}$.
\begin{definition}[Min--max sequences and critical pairs]
\label{def:minmax-families}
A sequence $\{\Omega_i\}\subset\Pi(M)$ is a \emph{min--max sequence} if
\[
\lim_{i\to\infty}\sup_{t\in[0,1]}\Aph(\Omega_i(t))=\Wph.
\]
Its \emph{critical pairs} are the pairs $(E,V)\in\Cacc\times\mathcal V(M)$ for which there exist $i_j\to\infty$ and $t_j\in[0,1]$ such that
\[
\Omega_{i_j}(t_j)\to E,
\qquad
|\partial\Omega_{i_j}(t_j)|\rightharpoonup V,
\qquad
\Aph(\Omega_{i_j}(t_j))\to\Wph.
\]
\end{definition}

\begin{theorem}[Anisotropic PMC min--max theorem]
\label{thm:apmc-minmax}
Let $(M^3,g)$ be a closed Riemannian manifold, let $\Phi$ be an elliptic integrand satisfying
\eqref{eq:elliptic-class}, and Let $h\in C^\infty(M)$ be $\Phi$-good and satisfy
$\int_M h\,\mathrm{d}V_g\geq 0$. Let $\{\Omega_i\}\subset\Pi(M)$ be a min--max sequence with one common
mass-ratio bound, i.e., there exist $C>0$ and $r_0>0$ such that
\[
\sup_{i\in\mathbb N}\sup_{t\in[0,1]}
\operatorname{MR}_{\Phi,r_0}(\Omega_i(t))
\le C
\]

Then there exists a critical pair $(\Omega,V)$ of $\{\Omega_i\}$ such that
\[
V=|\Sigma|,
\qquad
\Sigma=\partial\Omega,
\]
where $\Sigma$ is a nontrivial, smooth, closed, almost embedded
surface with prescribed anisotropic mean curvature $h$
with respect to the outward unit normal of $\Omega$.
Moreover,
\[
\mathcal A_{\Phi,h}(\Omega)=\mathbf W_{\Phi,h}.
\]
\end{theorem}

\vspace{15pt}
We prove Theorem~\ref{thm:apmc-minmax} at the end of this section, after establishing the existence of almost-minimizing varifolds and the regularity results.

\subsection{Stability, curvature estimates, and compactness}
\label{subsec:stable-estimates}

In this subsection, we establish the local curvature estimate and compactness for $h$-stable $h$-APMC boundaries.
\begin{theorem}[Stable curvature estimate]
\label{thm:stable-curvature}
Let $\Phi$ be an elliptic integrand on $M$ satisfying
\eqref{eq:elliptic-class},$U\subset M$ be open and let $\Sigma=\partial\Omega$ be a smooth almost embedded $h$-APMC boundary which is $h$-stable in $U$. If
\[
\mathcal H_g^2(\Sigma)\le D,
\qquad
\|h\|_{C^2(U)}\le H_2,
\]
then there exists $C=C(M,g,\lambda,H_2,D)$ such that
\begin{equation}
|A_\Sigma(x)|\le\frac{C}{\dist_g(x,\partial U)},
\qquad x\in\Sigma\cap U.
\label{eq:stable-curvature}
\end{equation}
\end{theorem}
\begin{proof}
The  argument in the proof of \cite[Theorem~2.8]{DD24}
gives uniform interior  mass-ratio bounds for the $h$-stable APMC boundaries. Indeed, this part of the argument uses only the stability inequality \eqref{eq:pmc-stability-prelim},it does not
require the anisotropic mean curvature to be constant or any additional consequence of $c$-stability.

Suppose that the conclusion were fail,then there exists a sequence
$\Sigma_j=\partial\Omega_j$, and points
$p_j\in\Sigma_j$ and scales $r _j=|A_{\Sigma_j}(p_j)|^{-1}\to 0$
such that
\[
r_j^{-1} \operatorname{dist}_g(p_j,\partial U) \to \infty.
\]

Now suppose after passing to a subsequence, $p_j\to p\in M$. The following is the standard blow-up argument as in
\cite[Lemma~1 and the proof of Theorem~1]{Whi87a}.
Identify $T_{p_j}M$ with $\mathbb R^3$ by orthonormal frames and rescale
using $\Psi_j(y)=\exp_{p_j}(r_jy)$. The rescaled surfaces have
second fundamental form of norm one at the origin and uniformly
bounded curvature on every fixed ball. The rescaled metrics converge
to the Euclidean metric, and the rescaled anisotropies converge to
the constant-coefficient integrand $\Phi(p,\cdot)$.

\vspace{5pt}
Thus after passage to a subsequence,the rescaled surface locally $C^2$ and graphically
convergence to a smooth closed  almost embedded stable $\Phi(p,\cdot)$-minimal surface $\Sigma_\infty\subset\mathbb R^3$. Strong maximum principle \cite{Ser70} implies $\Sigma_\infty$ is embedded and hence is two-sided. Moreover,
\[
 |A_{\Sigma_\infty}(0)|=1,
 \qquad \mathcal H^2(\Sigma_\infty\cap B_R(0))\leq C_*R^2
 \quad\text{for every }R>0,
\]
where $C_*$ is independent of $R$.White's Bernstein theorem
\cite[Theorem~4.1]{Whi91a} now implies that $\Sigma_\infty$ is a
plane, contradicting $|A_{\Sigma_\infty}(0)|=1$. This proves
\eqref{eq:stable-curvature}.
\end{proof}

The preceding curvature estimate gives the following compactness
corollary.
\begin{corollary}[Stable compactness]
\label{thm:stable-compactness}
Let $U\subset M$ be open, let $h_j\to h_\infty$ in $C^\infty(U)$, and let
$\Sigma_j=\partial E_j\subset U$ be smooth almost embedded $h_j$-stable APMC boundaries satisfying
\[
\sup_j\mathcal H_g^2(\Sigma_j)<\infty.
\]
After passing to a subsequence:
\begin{enumerate}[label=(\roman*)]
\item if $h_\infty$ is $\Phi$-good, then $\Sigma_j$ converges locally smoothly in $U$ to a smooth almost embedded $h_\infty$-stable APMC boundary $\Sigma$; moreover
\[
\Theta(\Sigma,x)=1\quad(x\in\mathcal R(\Sigma)),
\qquad
\Theta(\Sigma,x)=2\quad(x\in\mathcal S(\Sigma));
\]
\item if $h_\infty\equiv0$, then $\Sigma_j$ converges locally smoothly in $U$, with integer multiplicity, to a smooth embedded $\Phi$-stable surface $\Sigma$.
\end{enumerate}
\end{corollary}
\begin{proof}
If $h$ is $\Phi$-good, then $\Sigma$ has no $\Phi$-minimal component. Then the proof follows exactly as in \cite[Theorem~3.19]{ZZ20}, where we use Theorem~\ref{thm:stable-curvature} in place of \cite[Theorem~3.6]{ZZ20}.

\end{proof}

\subsection{almost-minimizing varifolds and replacement  with controlled mass ratio}

In this subsection, we introduce the notion of $h$-almost minimizing varifolds and develop the corresponding replacement theory with controlled mass ratio.

\begin{definition}[$(\eps,\delta)$-almost minimizing sets and $h$-almost minimizing varifolds]
Let $U\subset M$ be open and $\eps,\delta>0$. Define
$\mathfrak a_{\Phi,h}(U;\eps,\delta)$ as the set of all $E\in\Cacc$ such that every finite sequence $E_0,\ldots,E_m\in\Cacc$ satisfying
\begin{enumerate}[label=(\roman*)]
\item $E_0=E$;
\item $E_i\triangle E\subset U$ modulo null sets, $i=1,\ldots,m$;
\item $\Vol(E_i\triangle E_{i-1})\le\delta$, $i=1,\ldots,m$;
\item $\Aph(E_i)\le\Aph(E)+\delta$, $i=1,\ldots,m$,
\end{enumerate}
satisfies
\[
\Aph(E_m)\ge\Aph(E)-\eps.
\]
A varifold $V\in\mathcal V(M)$ is \emph{$h$-almost minimizing for $\Phi$ in $U$} if there exist
\[
\eps_j\downarrow0,
\qquad
\delta_j\downarrow0,
\qquad
E_j\in\mathfrak a_{\Phi,h}(U;\eps_j,\delta_j),
\qquad
\bfF_V(|\partial E_j|,V)\le\eps_j.
\]
When the sequence $E_j$ is specified, we say that $V$ is $h$-almost minimizing in $U$ \emph{by means of $E_j$}.
\end{definition}
The following notion  was first introduced in
\cite{Li23} in the study of Almgren--Pitts width realizations.
\begin{definition}
For $p\in M$ and $0<s<r<\operatorname{inj}(M)$, set
\[
\an(p;s,r):=B_r(p)\setminus\overline{B_s(p)}.
\]
A collection
\[
\mathscr C=\{\an(p;s_1,r_1),\ldots,\an(p;s_L,r_L)\}
\]
is $L$-admissible if
\[
2r_{j+1}<s_j,
\qquad j=1,\ldots,L-1.
\]
A varifold $V$ is \emph{$h$-almost minimizing in $\mathscr C$} (by means of $E_j$) if it is $h$-almost minimizing in at least one $A\in\mathscr C$ (by means of $E_j$). 
\end{definition}

\begin{definition}
We say a $V$ is
$h$-almost minimizing in small annuli (by means of
$\{E_j\}$) if, for every $p\in M$, there exists
$r_{\mathrm{am}}(p)>0$ such that $V$ is $h$-almost
minimizing in $\operatorname{An}(p;s,r)$
(by means of $\{E_j\}$) whenever
$0<s<r\le r_{\mathrm{am}}(p)$.
\end{definition}

The following lemma is proved by the same
argument as in
\cite[Lemma~4.5]{DDL24}, and we omit its proof.

\begin{lemma}
Let $L$ be a positive integer. If a varifold $V$ is
$h$-almost minimizing in every $L$-admissible
collection of annuli (by means of $\{E_j\}$),
then $V$ is $h$-almost minimizing in small annuli
(by means of $\{E_j\}$).
\end{lemma}

We now apply the Pitts combinatorial argument \cite{Pit81}to obtain an almost-minimizing critical varifold.
\begin{theorem}[Existence of Almost-minimizing varifolds]
\label{thm:AM-extraction-417}
Let $\{\Omega_i\}\subset\Pi(M)$ be a min--max sequence with one common
mass-ratio bound, i.e., there exist $C>0$ and $r_0>0$ such that
\[
\sup_{i\in\mathbb N}\sup_{t\in[0,1]}
\operatorname{MR}_{\Phi,r_0}(\Omega_i(t))
\le C
\]Then there exists a  critical pair $(E,V)$ such that
\begin{enumerate}[label=(\roman*)]
\item $V$ has the same mass-ratio bound as $\{\Omega_i\}$ ;
\item $V$ has $c$-bounded $\Phi$-first variation,where $c=\|h\|_{C^0}$;
\item for some universal  constant $\bar L\le27$, $V$ is $h$-almost minimizing in every $\bar L$-admissible collection of annuli, by means of $\Omega_{i_j}(t_j)$ for some $i_j\to\infty$ and $t_j\in[0,1]$ such that
\[
\Aph(\Omega_{i_j}(t_j))\to\Wph.
\].
\end{enumerate}

\end{theorem}
\begin{proof}
The existence of an almost-minimizing critical varifold is obtained
from Pitts' combinatorial argument \cite[Theorem~4.10]{Pit81}.
The contradiction is carried out directly on slices of the given
min--max sequence whose $\mathcal A_{\Phi,h}$-values approach
$\mathbf W_{\Phi,h}$. If no critical varifold arising as a limit of
such slices were $h$-almost minimizing in every $L$-admissible
collection of annuli, then the failure of the $h$-almost minimizing
condition would provide, for each such slice, a deformation supported
in a suitable annulus that decreases $\mathcal A_{\Phi,h}$ by a
definite amount. Pitts' combinatorial selection combines these
deformations while preserving the sweepout class, producing a
competing sequence whose maximal $\mathcal A_{\Phi,h}$-value is
strictly below $\mathbf W_{\Phi,h}$, a contradiction. In Pitts'
construction one may take
\[
L=(3^m)^{3^m},
\]
where $m$ is the dimension of the parameter space. Since our sweepouts
are one-parameter, $m=1$, and hence $L=27$.Thus, (iii) is proved, and (i) follows directly from varifold convergence.

Now suppose $V$ does not satisfy (ii), after using a nonnegative partition of unity, there exists a
vector field $Y$ supported in some $B_{r_{0}}(p)$ for which
\[
 \delta_\Phi V(Y)+c\int|Y|\,\dd\|V\|<0.
\]
Let $\eta_\rho=1$ on $B_\rho(p)$, vanish outside
$B_{2\rho}(p)$, and satisfy $|D\eta_\rho|\leq C/\rho$.
The  mass ratio bound together with \cite[(2.3)]{DDL24}
give
\begin{equation}
 |\delta_\Phi V(\eta_\rho Y)|
       +c\int|\eta_\rho Y|\,\dd\|V\|
 \leq C\bigl(\rho^2\|DY\|_\infty+
              \rho\|Y\|_\infty+H\rho^2\|Y\|_\infty\bigr)
 \longrightarrow0.
 \label{eq:puncture-first-variation}
\end{equation}
Choose $\rho>0$ sufficiently small and set $Z:=(1-\eta_\rho)Y$. Then
$\operatorname{spt}(Z)\subset A(p,\rho,r_0)$ and, for some $\gamma>0$,
\[
\delta_\Phi V(Z)
    +H\int |Z|\,d\|V\|
    <-2\gamma .
\]
Let $\{\varphi_t\}_{t\in\mathbb R}$ be the flow generated by $Z$.
By continuity, there exist $\alpha>0$ and $\varepsilon>0$ such that,
whenever $\mathbf F(|\partial\Omega|,V)<\varepsilon$,
\[
\frac{d}{dt}
\mathcal A_{\Phi,h}\bigl(\varphi_t(\Omega)\bigr)
<-\gamma ,
\qquad 0\le t\le\alpha .
\]
Hence, for $0\le t_1<t_2\le\alpha$,
\[
\mathcal A_{\Phi,h}\bigl(\varphi_{t_2}(\Omega)\bigr)
<
\mathcal A_{\Phi,h}\bigl(\varphi_{t_1}(\Omega)\bigr),
\qquad
\mathcal A_{\Phi,h}(\Omega)
-
\mathcal A_{\Phi,h}\bigl(\varphi_\alpha(\Omega)\bigr)
\ge \gamma\alpha .
\]

Write $\Omega_j:=\Omega_{i_j}(t_j)$. Since
$|\partial\Omega_j|\to V$ and $\Omega_j$ is
$(\varepsilon_j,\delta_j)$-$h$-almost minimizing in every
$\bar L$-admissible collection of annuli, arguing as in
\cite[Lemma~4.5]{DDL24}, and choosing $r_0>0$ sufficiently small if
necessary, we may assume that $\Omega_j$ is
$(\varepsilon_j,\delta_j)$-$h$-almost minimizing in
$A(p,\rho,r_0)$. For all sufficiently large $j$,
$\mathbf F(|\partial\Omega_j|,V)<\varepsilon$ and
$\varepsilon_j<\gamma\alpha$.
After subdividing $[0,\alpha]$ sufficiently finely, the sequence
\[
\Omega_j,\ 
\varphi_{t_1}(\Omega_j),\ldots,
\varphi_{t_N}(\Omega_j),
\qquad 0=t_0<t_1<\cdots<t_N=\alpha ,
\]
satisfies the admissible deformation conditions in
$A(p,\rho,r_0)$, while
\[
\mathcal A_{\Phi,h}\bigl(\varphi_{t_N}(\Omega_j)\bigr)
<
\mathcal A_{\Phi,h}(\Omega_j)-\varepsilon_j .
\]
Thus
$\Omega_j\notin
\mathfrak a_{\Phi,h}
(A(p,\rho,r_0);\varepsilon_j,\delta_j)$,
a contradiction. Therefore $V$ satisfies~(ii).
\end{proof}

We next recall the following replacement notation, which will play an important role in the regularity argument below.
\begin{definition}[$h$-replacement]
Let $U\subset M$ be open, $K\Subset U$, and let $V$ be $h$-almost minimizing in $U$. A varifold $V^*\in\mathcal V(M)$ is an \emph{$h$-replacement for $V$ in $K$} if 
\begin{enumerate}[label=(\roman*)]
\item
$V^*\llcorner(M\setminus K)=V\llcorner(M\setminus K)$;
\item
$|\Fphi(V^*)-\Fphi(V)|\le c\Vol(K)$,where $c=\|h\|_{C^0}$;
\item 
$V^*=|\Sigma^*|$ in $\operatorname{int}K$,
where $\Sigma^*$ is a smooth almost embedded $h-$ stable APMC boundary .

\end{enumerate}
\end{definition}
\begin{lemma}
\label{lem:replacement-minimization}
Let $U\subset M$ be open, $K\subset U$, $\varepsilon,\delta>0$, $\|h\|_{C^0(M)}\le\Lambda$,and
$E\in\mathfrak a_{\Phi,h}(U;\varepsilon,\delta)$. Define
$\mathcal C_E^\delta(K)$ to be the collection of all
$\widetilde E\in\mathcal C(M)$ for which there exists a finite sequence
$E=E_0,E_1,\ldots,E_m=\widetilde E$
such that
\begin{enumerate}[label=(\alph*)]
    \item $E_i\triangle E\subset K$ for every $i=1,\ldots,m$;
    \item $\Vol(E_i\triangle E_{i-1})\leq\delta$ for every
    $i=1,\ldots,m$;
    \item $\mathcal A_{\Phi,h}(E_i)
    \leq\mathcal A_{\Phi,h}(E)+\delta$ for every
    $i=1,\ldots,m$.
\end{enumerate}
Then there exists $E^*\in\mathcal C_E^\delta(K)$ such that
\begin{enumerate}[label=(\roman*)]
    \item $E^*$ locally minimizes $\mathcal A_{\Phi,h}$ in
    $\operatorname{int}K$;
    \item $E^*\in\mathfrak a_{\Phi,h}(U;\varepsilon,\delta)$;
    \item if $\operatorname{MR}_{\Phi,r_0}(E)\leq C_0$, then
    \[
    \operatorname{MR}_{\Phi,r_0}(E^*)\leq C
    \]
    for a constant
    $C=C(C_0,M,g,\lambda,\Lambda)$;
    \item
    \[
    \mathcal A_{\Phi,h}(E)-\varepsilon
    \leq
    \mathcal A_{\Phi,h}(E^*)
    \leq
    \mathcal A_{\Phi,h}(E)+\delta.
    \]
\end{enumerate}
\end{lemma}
\begin{proof}
The proof of (i)(ii)(iii) is the same as that of \cite[Lemma 4.8]{DDL24}, with
$\Phi^c$ replaced by $\mathcal A_{\Phi,h}$, every term of the form
$c\,\Vol(E)$ replaced by $\int_E h\,\dd V_g$, and the constant $c$ in
the estimates for inequalities replaced by $H:=\|h\|_{C^0}$. All the arguments and estimates in
\cite[Lemma 4.8]{DDL24} carry over with these modifications, and the
conclusions follow.
Finally, since $E^*\in\mathcal C_E^\delta(K)$ and $K\subset U$, the
upper bound in (iv) follows from condition~(c) in the definition of
$\mathcal C_E^\delta(K)$, while the lower bound follows directly from
$E\in\mathfrak a_{\Phi,h}(U;\varepsilon,\delta)$.
\end{proof}

We conclude this section by establishing the existence of replacements for
$h$-almost minimizing varifolds and recording some of their basic properties.

\begin{proposition}[Existence and properties of $h$-replacements]
\label{prop:replacement-properties}
Let $h$ be $\Phi$-good, and $\|h\|_{C^0(M)}\le\Lambda$.let $U\subset M$ be open and $K\Subset U$, and suppose that $V$ is $h$-almost minimizing in $U$ by means of
\[
E_j\in\mathfrak a_{\Phi,h}(U;\eps_j,\delta_j),
\qquad
\eps_j,\delta_j\downarrow0,
\]
with  common mass-ratio bound $C_0$ for the $E_j$. Then there exists $C=C(C_0,M,g,\lambda,\Lambda)$ and a $h$-replacement $V^*$ for $V$ in $K$ such that
\begin{enumerate}[label=(\roman*)]
\item 
$V^*$ is $h$-almost minimizing for $\Phi$ in $U$ by means of
$E_j^*$, $\varepsilon_j$, and $\delta_j$, for some
$E_j^*\in\mathcal C(M)$ such that $E_j^*$ is $h$-stable in
$\operatorname{int}K$, locally minimizes $\mathcal A_{\Phi,h}$ in
$\operatorname{int}K$, and has mass-ratio bound $C$.
In particular, $V^*$ also has the same mass-ratio bound $C$.
\item If $V$ has $c$-bounded anisotropic first variation,then so does $V^*$.

\end{enumerate}

\end{proposition}
\begin{proof}
By the regularity results of \cite{Bom82,DM15}, each locally $\mathcal A_{\Phi,h}$-minimizing set provided by
Lemma~\ref{lem:replacement-minimization} is equivalent to an open set whose reduced boundary is a
smooth $h$-stable APMC surface in $\operatorname{int}K$.
With this regularity result, the proof is the same as that of \cite[Proposition~5.8]{ZZ19},
with \cite[Lemma~5.7]{ZZ19} replaced by Lemma~\ref{lem:replacement-minimization} and using
Corollary~\ref{thm:stable-compactness} for the compactness of the $h$-stable APMC boundaries.
\end{proof}

\subsection{Regularity of almost-minimizing varifolds}

\begin{theorem}[Regularity]
\label{thm:APMC-regularity}
Let $\Lambda>0$, let $h\in C^\infty(M)$ be a
$\Phi$-good function with $\|h\|_{C^0(M)}\le\Lambda$,
and let $V$ be the varifold constructed in
Theorem~\ref{thm:AM-extraction-417}. Then $V=|\Sigma|$, where $\Sigma$
is a smooth almost embedded $h$-APMC surface.
Moreover,
\[
\Theta(V,x)=1 \quad\text{for every }x\in R(\Sigma),
\qquad
\Theta(V,x)=2 \quad\text{for every }x\in S(\Sigma).
\]
Furthermore, $\Sigma$ has has a mass-ratio bound $C$ below scale $r_0$,where
$C=C(M,g,\lambda,\Lambda)$.

\end{theorem}

\begin{proof}
We follow the proof strategy of
\cite[Theorem~5.1]{DDL24}, with the necessary
modifications for the prescribed anisotropic
mean curvature setting.

Let $V$ be the $2$-varifold obtained in
Theorem~\ref{thm:AM-extraction-417}.
We establish its regularity in a neighborhood of each
$x\in\operatorname{spt}\|V\|$.
Fix such an $x$ and choose $\rho>0$ sufficiently small
that $2\rho<r_{\mathrm{am}}(x)$.
By the $h$-almost minimizing property, $V$ is
$h$-almost minimizing in every annulus
$\operatorname{An}(x;s,r)$ with
$0<s<r\le 2\rho$.
In particular, the corresponding replacements are
available by Proposition~\ref{prop:replacement-properties}.

Let $V'$ be an $h$-replacement of $V$ in
$\operatorname{An}(x;\rho,2\rho)$, whose restriction
to this annulus is induced by a smooth, almost
embedded, $h$-stable surface $\Sigma'$.
We next select a suitable radius for the second
replacement. By Theorem~\ref{thm:touching},
the touching set $S(\Sigma')$ has locally finite
one-dimensional Hausdorff measure.
The coarea formula applied to
$r_x(y)=d_g(x,y)$ yields
\[
\begin{aligned}
\int_a^b
\mathcal H^0\bigl(
S(\Sigma')\cap\partial B_t(x)
\bigr)\,dt
&\le
\mathcal H_g^1\bigl(
S(\Sigma')\cap\operatorname{An}(x;a,b)
\bigr)
<\infty
\end{aligned}
\]
whenever $\rho<a<b<2\rho$.
Together with Sard's theorem, this allows us to
choose $t\in(\rho,2\rho)$ such that
$\partial B_t(x)$ is transverse to  $\Sigma'$ and
\begin{equation}\label{eq:count touching point}
\#\bigl(
\mathcal{S}(\Sigma')\cap\partial B_t(x)
\bigr)<\infty.
\end{equation}
In fact, both conditions are satisfied for almost
every $t\in(\rho,2\rho)$.

Fix any $s\in(0,\rho)$ and take an $h$-replacement
$V''$ of $V'$ in $\operatorname{An}(x;s,t)$.
Its restriction to this annulus is represented by
a smooth, almost embedded, $h$-stable surface
$\Sigma''$.
Moreover, both $V'$ and $V''$ are integral
$2$-varifolds by \cite[Lemma~4.15]{DD24}.

\medskip
\noindent\textbf{Step 1.}
We claim that for every
$y\in \Sigma'\cap\partial B_t(x)$,
there exists a sufficiently small radius $r>0$
such that
\[
\Sigma''\cap B_t(x)\cap B_r(y)
=
\Sigma'\cap B_t(x)\cap B_r(y).
\]
Since the background can be absorbed into the anisotropic, we will
assume without loss of generality that the
ambient space is $\mathbb R^3$.

\medskip
\noindent\textit{Case 1: $y\in \mathcal{R}(\Sigma')$.}
Set $H_0:=\|h\|_{C^0(M)}$ and fix $0<\varepsilon\ll1$.
Let $r>0$ small enough such that
\[
\sup_{q\in\Sigma'\cap B_{2r}(y)}
\operatorname{dist}_{G}(T_q\Sigma',T_y\Sigma')<\varepsilon,
\qquad B_{2r}(y)\cap \mathcal{S}(\Sigma')=\emptyset.
\]
Choose
$z\in\Sigma'\cap B_r(y)\setminus\overline{B_t(x)}$.
As in the  construction of
\cite[Theorem~5.1, Step~1, Case~1]{DDL24}, choose a convex
set $\widetilde{\mathcal C}$ bounded by two spherical caps
$\Gamma_+$ and $\Gamma_-$ meeting along a smooth circle
$S\ni z$ at angle $3\varepsilon$ and furthermore the following hold:
\[
T_zS=T_z(\Sigma'\cap\partial\widetilde{\mathcal C}),
\quad
\Sigma'\cap\bigl(B_r(y)\setminus
\operatorname{An}(x;s,t)\bigr)\subset\widetilde{\mathcal C}.
\]
Define 
\[
\mathcal C:=\widetilde{\mathcal C}\setminus
\left(y+\tfrac12(\widetilde{\mathcal C}-y)\right),\qquad
K:=\overline{\mathcal C}
\]
By taking $r$ smaller, if necessary, we can find
$0<r_1<r_2<r_{\mathrm{am}}(y)$ with
\[
K\subset \operatorname{An}(y;r_1,r_2).
\]
By Proposition~\ref{prop:replacement-properties} we have $V''$
$h$-almost minimizing in $\operatorname{An}(y;r_1,r_2)$. Hence it admits a replacement
$V'''$ in $K$, satisfying
\[
V'''\llcorner(M\setminus K)=V''\llcorner(M\setminus K),
\qquad
\|\delta_\Phi V'''\|\le H_0\|V'''\|.
\]
In $\mathcal C$, $V'''$ is represented by a smooth,
$h$-stable, almost embedded APMC boundary.
For any varifold $Q$, write $\operatorname{TanVar}(q,Q)$ for the set of
weak limits of $(\eta_{q,r_j})_\#Q$ as $r_j\downarrow0$,
where $\eta_{q,r}(p)=(p-q)/r$.
The mass-ratio bound then yield
$\operatorname{TanVar}(z,V''')\ne\varnothing$.
Up to rotation,we can assume that
\[
\begin{gathered}
T_zK=\mathcal W_\varepsilon
:=\{p\in\mathbb R^3:p_1\ge0,
\ |p_3|\le\tan(3\varepsilon)p_1\},\\
\ell:=T_zS=\operatorname{span}(e_2),\qquad
P:=T_z\Sigma'=\operatorname{span}(e_1,e_2).
\end{gathered}
\]
Since $V''=|\Sigma'|$ in a neighborhood of $z$
and $V'''=V''$ outside $K$, every
$W\in\operatorname{TanVar}(z,V''')$ satisfies
\begin{equation}\label{eq:pmc-gluing-exterior-halfplane}
W\llcorner\{p_1\leq0\}=|P|\llcorner\{p_1\leq0\}.
\end{equation}
We next examine the structure
of tangent varifolds in the half-space
$\{p_1>0\}$ through the following lemma.
\begin{lemma}
\label{lem:pmc-extremal-wedge}
Let $H:=\{p_1>0\}$ and
$\Phi_z(\nu):=\Phi(z,\nu)$.
For some $W\in\operatorname{TanVar}(z,V''')$ there exist
$\nu_1,\nu_2\in\mathbb S^2\cap\operatorname{span}(e_1,e_3)$
such that
\[
\operatorname{spt}\|W\|\cap H\subset
H\cap\{\langle\nu_1,p\rangle\ge0\}
\cap\{\langle\nu_2,p\rangle\le0\}
\subset\mathcal W_\varepsilon.
\]
Moreover, writing $Q_i:=\nu_i^\perp$, either
\begin{equation}\label{eq:pmc-extremal-two}
\nu_1\ne\nu_2,\qquad
|Q_1|\llcorner H+|Q_2|\llcorner H
\le W\llcorner H,
\end{equation}
or
\begin{equation}\label{eq:pmc-extremal-one}
\nu_1=\nu_2,\qquad
W\llcorner H=|Q_1|\llcorner H.
\end{equation}
\end{lemma}

\begin{proof}
Since $V'''=V''$ outside $K$ and
$V''=|\Sigma'|$ near $z$, we have
$
\operatorname{spt}\|V'''\|\cap U
\subset (K\cap U)\cup(\Sigma'\cap U)
$
for some neighborhood $U$ of $z$.
Consequently, for every
$W\in\operatorname{TanVar}(z,V''')$,we have
$
\operatorname{spt}\|W\|
\subset T_zK\cup P.
$
Since $P\cap\{p_1\ge0\}\subset T_zK$, it follows that
$
\operatorname{spt}\|W\|\cap\{p_1\ge0\}
\subset T_zK=\mathcal W_\varepsilon.
$
Then this fact supplies the
counterpart of \cite[(5.12)--(5.13)]{DDL24}.
The entire argument following (5.13) in the proof of
\cite[Lemma~5.2]{DDL24} carries over verbatim to our setting,
as those argument only involves stationary varifolds
with respect to the frozen integrand.
\end{proof}

\begin{lemma}[Multiplicity-one planar tangent]
\label{lem:pmc-planar-tangent}
There exists $\varepsilon_0=\varepsilon_0(\lambda)>0$
such that, if $\varepsilon<\varepsilon_0$, the varifold
$W$ in Lemma~\ref{lem:pmc-extremal-wedge}
satisfies
\[
W=|P|=|T_z\Sigma'|.
\]
\end{lemma}
\begin{proof}
Replace \cite[Lemma~5.2]{DDL24} with 
Lemma~\ref{lem:pmc-extremal-wedge} in the proof of
\cite[Lemma~5.3]{DDL24},then the entire argument in the proof of
\cite[Lemma~5.3]{DDL24} carries over verbatim to our setting,
as those argument only involves stationary varifolds
with respect to the frozen integrand.
\end{proof}

By Lemma~\ref{lem:pmc-planar-tangent}, there is a sequence
$\tau_j\downarrow0$ such that
\[
 W_j:=(\eta_{z,\tau_j})_\#V'''\rightharpoonup |P|,
 \qquad P:=T_z\Sigma'.
\]
 Write
\[
 \Phi_j(y,\nu):=\Phi(z+\tau_jy,\nu),
 \qquad H_0:=\|h\|_{C^0(M)}.
\]
Then the integrands $\Phi_j$ converge in $C^3$ on compact sets to the frozen
integrand $\Phi_z(\nu)=\Phi(z,\nu)$.Thus there exists a constant $\lambda\geq 1$ such that $\Phi_j$ satisfying the quantitative in
\eqref{eq:elliptic-class}.We show that, for all sufficiently large $j$, a restriction of $W_j$ on a neighborhood of $0$
satisfies the hypotheses of Theorem~\ref{thm:allard-uniform-mass-lower}.

\vspace{10pt}
\emph{Condition (v).}
The following proof proceeds as in
\cite[Lemma~4.15]{DD24}.
The key point is that, by
Proposition~\ref{prop:replacement-properties},
there exists $r>0$ such that $V'''$ remains
$h$-almost minimizing in $B_{4r}(z)$.
 We claim that there are
$\kappa_0>0$,such that
\begin{equation}\label{eq:local-uniform-lower-mass-at-z}
 \|V'''\|(B_\rho(q))\ge\kappa_0\rho^2
 \quad\text{for all }q\in\operatorname{spt}\|V'''\|\cap B_{3r}(z),
 \quad 0<\rho\le r.
\end{equation}
Indeed, given such $q$ and $\rho$, 
by Proposition~\ref{prop:replacement-properties},we can take an $h$-replacement
$\widehat V_{q,\rho}$ of $V'''$ in $\overline{B_{3\rho/4}(q)\setminus B_{\rho/4}(q)}$. Since
$\widehat V_{q,\rho}=V'''$ on $B_{\rho/4}(q)$ and
$q\in\operatorname{spt}\|V'''\|$,we claim that $\spt\widehat V_{q,\rho}\bigcap \partial B_{\rho/2}(q)\neq \emptyset$.Assume the claim fails, there exists some
$
r_1<\frac{\rho}{2}
$
such that $\spt\widehat V_{q,\rho}$ touches $\partial B_{r_1}(y)$ from the inside. More
precisely,for some sufficiently small $\varepsilon>0$,
\[
\spt\widehat V_{q,\rho}\cap B_{\rho}(q)\cap\partial B_{r_1}(q)\neq\emptyset,\quad \spt\widehat V_{q,\rho}\cap B_{\rho}(q)\cap \bigl[\overline{B_{r_1+\varepsilon}(y)}\setminus\overline{B_{r_1(y)}}\bigr]=\emptyset.
\]
But this conclusion,by \cite[Lemma~4.5]{DRKS20},
contradicts the maximum principle in
\cite[Proposition~2(iii)]{DDH19},thus there exists $x\in\spt\widehat V_{q,\rho}\bigcap \partial B_{\rho/2}(q)$, the uniform mass ratio lower bound for $h-$ stable APMC surface (see \cite[p.~305, Theorem, (9)]{All83})
gives a constant $a_0>0$, independent of $q$ and $\rho$, for which
\[
 \|\widehat V_{q,\rho}\|(B_\rho(q))\ge\|\widehat V_{q,\rho}\|(B_\rho/4(x))\ge a_0\rho^2.
\]

\noindent By the definition of replacement we have
$\widehat V_{q,\rho}=V'''$ outside $K_{q,\rho}\Subset B_\rho(q)$,
\[
 \mathbf F_\Phi(V''';B_\rho(q))
 \ge \mathbf F_\Phi(\widehat V_{q,\rho};B_\rho(q))
      -H_0\operatorname{Vol}(K_{q,\rho}).
\]
Using $\lambda^{-1}\le\Phi\le\lambda$ and
$\operatorname{Vol}(K_{q,\rho})\le C\rho^3$, we obtain
\[
 \|V'''\|(B_\rho(q))
 \ge\lambda^{-2}a_0\rho^2-C\lambda^{-1}H_0\rho^3.
\]
By choosing $r>0$ sufficiently small, the $O(\rho^3)$ term can be absorbed into the leading-order $\rho^2$ term. This completes the proof of \eqref{eq:local-uniform-lower-mass-at-z}.

Rescaling \eqref{eq:local-uniform-lower-mass-at-z}, we obtain, for all sufficiently large $j$,

\begin{equation}\label{eq:rescaled-uniform-lower-mass-at-z}
 \|W_j\|(B_t(a))
 =\tau_j^{-2}\|V'''\|(B_{\tau_jt}(z+\tau_ja))
 \ge\kappa_0t^2
\end{equation}
for every $a\in\operatorname{spt}\|W_j\|\cap B_3(0)$ and
$0<t\le1$. In particular, this implies condition (v) hold for all sufficiently large $j$.

\emph{Conditions (ii)--(iv).}
Set $\Omega=B_2(0)$, and fix $s=1$ and $\alpha\in(0,1)$.
Let $r_*=r_*(n,\alpha,C,\lambda)>0$ be the constant
in Theorem~\ref{thm:allard-uniform-mass-lower}, where $C$ denotes the uniform mass-ratio
upper bound for $V'''$.
Fix $0<r_0<\min\{r_*,1/8\}$ and take $r=r_0$ in conditions (ii)--(iv). Then, for all sufficiently large $j$, these conditions follow directly from the uniform mass-ratio upper bound
for $V'''$, the varifold convergence $W_j\rightharpoonup|P|$,
and the fact $\|\delta_\Phi W_j\|\le \tau_jH_0\|W_j\|$ .

\emph{Condition (i).}
The lower bound \eqref{eq:rescaled-uniform-lower-mass-at-z} and
$W_j\rightharpoonup|P|$ in fact imply
\[
 \sup_{a\in\operatorname{spt}\|W_j\|\cap B_2(0)}
       \operatorname{dist}(a,P)\longrightarrow0.
\]
Otherwise, a sequence of point ${x_j} \in \operatorname{spt}\|W_j\|\cap B_2(0)$ remaining a fixed positive distance from
$P$ would carry, by \eqref{eq:rescaled-uniform-lower-mass-at-z} and weak
convergence $W_j\rightharpoonup |P|$, a
fixed positive mass of $|P|$ in a ball disjoint from $P$. 
Writing points in $\mathbb R^3$ as
$(s,t)\in\mathbb R^2\times\mathbb R$, with $t$
denoting the last coordinate, then
\[
 B_{4r_0}^2(0)\times[-1,1]\subset\Omega,
 \qquad
 \operatorname{spt}\|W_j\|\cap\mathcal C_{4r_0}
 \subset\{|t|\le1\}
\]
for all sufficiently large $j$, where
$\mathcal C_{4r_0}=\Omega\cap(B_{4r_0}^2(0)\times\mathbb R)$.
Thus, condition (i) holds for $W_j\llcorner\Omega$ for all sufficiently large $j$.

\emph{Condition (vi).}
For this fixed $r_0$ above, the weak convergence $W_j\rightharpoonup|P|$, together with the
height bound just proved, gives
\[
 \mu_j^2:=(4r_0)^{-4}
       \int_{\mathcal C_{4r_0}}t^2\,d\|W_j\|\longrightarrow0.
\]
Moreover $4r_0\tau_jH_0\to0$. Consequently,
condition (vi) holds for $W_j\llcorner\Omega$ for all sufficiently large $j$.

\vspace{10pt}

We have verified all  hypotheses of
Theorem~\ref{thm:allard-uniform-mass-lower} for one fixed,
sufficiently large index $j_0$. It follows that
$\operatorname{spt}\|W_{j_0}\|$ is a single $C^{1,\alpha}$ graph in
$\mathcal C_{2r_0}$ and $\operatorname{spt}\|V'''\|$ is a
single $C^{1,\alpha}$ graph $G_z$ in a neighborhood $O_z$ of $z$.Then standard elliptic regularity gives $G_z$ is smooth and $z$ is a regular
point of $V'''$.By the fact
\[
V'''\llcorner O_z=|G_z|,
\qquad G_z\cap(O_z\setminus K)
=\Sigma'\cap(O_z\setminus K).
\]
we have $V'$ and $V'''$ glue smoothly at z,unique continuation then yields

\[
V'''\llcorner\mathcal C
=
|\Sigma'|\llcorner\mathcal C
\quad\text{on the connected component of }
\operatorname{spt}\|V'''\|
\text{ passing through }z.
\]

\noindent Shrinking $r$ further if necessary,we can assume that the two smooth caps of
$\partial\widetilde{\mathcal C}\setminus S$ satisfy
\[
\inf_{\partial\widetilde{\mathcal C}\setminus S}
H_\Phi[\partial\widetilde{\mathcal C},
\nu_{\rm out}.
]
\ge\frac{c_\varepsilon}{r}>2H_0.
\]

 \noindent where $\nu_{\rm out}$ denotes the outward unit normal
to $\widetilde{\mathcal C}$ on each spherical cap.

\vspace{5pt}
\noindent we claim that:
\begin{equation}\label{eq:pmc-lens-outside-approach}
\Sigma'\cap\partial\widetilde{\mathcal C}
\subset\operatorname{spt}\|V'''\|
\cap\partial\widetilde{\mathcal C}
\subset{\operatorname{spt}\|V'''\|
\cap(\mathbb R^3\setminus\widetilde{\mathcal C})}\subset\operatorname{spt}\|V''\|.
\end{equation}

\noindent Otherwise, $\operatorname{spt}\|V'''\|\cap\widetilde C$
would touch $\partial\widetilde C$ from the inside,
which, by \cite[Lemma~4.5]{DRKS20},
contradicts the maximum principle in
\cite[Proposition~2(iii)]{DDH19}.
.

\vspace{8pt}
Choose a nested family of lenses
$\{\widetilde{\mathcal C}_\tau\}_{\tau\in[0,1]}$
along a curve on  $\Sigma'$ joining $y$ to $z$,
with $\widetilde{\mathcal C}_0=\{y\}$ and
\[
\widetilde{\mathcal C}_1\subset
\bigcup_{\tau\in[0,1]}
\partial\widetilde{\mathcal C}_\tau.
\]
Repeating \eqref{eq:pmc-lens-outside-approach} along this
family yields
\[
\Sigma'\cap\widetilde{\mathcal C}_1
\subset\operatorname{spt}\|V''\|.
\]
Since both $V'$ and $V''$ are integral
$2$-varifolds by \cite[Lemma~4.15]{DD24},
we obtain in a neighborhood of $y$
\begin{equation}\label{eq:pmc-inner-varifold-dominates}
V''\ge|\Sigma'|.
\end{equation}
Thus for every $Z\in\operatorname{TanVar}(y,V'')$, writing
$\Phi_y(\nu)=\Phi(y,\nu)$, we infer
\[
R_Z:=Z-|T_y\Sigma'|\ge0,\qquad
\delta_{\Phi_y}R_Z=0,\qquad
\operatorname{spt}\|R_Z\|\subset\{p_1\ge0\}.
\]

The following lemma is a PMC analogue of the estimate in
\cite[proof of Proposition~5.4, (5.4)--(5.6)]{DD24} which establishes the uniqueness of the tangent varifolds  of $V''$.

\begin{lemma}
\label{lem:pmc-interface-wedge}
Under the preceding assumptions,there exists
a constant $L<\infty$ such that 

\begin{equation}\label{eq:pmc-spherical-wedge}
\begin{gathered}
\operatorname{spt}\|Z\|\cap\{p_1\ge0\}\subset\mathcal W_L
:=\{p_1\ge0,\ |p_3|\le Lp_1\},\\\qquad \text{and} \qquad R_Z=0
\end{gathered}
\end{equation}

\end{lemma}

\begin{proof}
After choosing suitable local coordinates
$(p_1,p_2,p_3)$ on a neighborhood $U$ of $y$,
centered at $y$, we may assume that
 ,we can assume that 
\[
\begin{gathered}
T_y\partial B_t(x)
=\{p_1=0\},\\
P:=T_y\Sigma'=\{p_3=0\},\\
\ell:=P\cap\{p_1=0\}=\{p_1=p_3=0\},
\end{gathered}
\]

\noindent Furthermore, write
$q=(p_2,p_3)$,there exits a smooth function $\psi$ on $span\{p_2,p_3\}$ such that
\[
\begin{gathered}
\partial B_t(x)\cap U=\{p_1=\psi(q)\},\\
B\cap U=\{p_1>\psi(q)\},\\
\psi(0)=0,\qquad D\psi(0)=0.
\end{gathered}
\]
Without of loss of generality,we can assume that $t$ small enough such that
\[
\kappa:=
\inf_{\partial B_t(x)\cap U}
H_\Phi[\partial B_t(x),\nu_{\rm out}]>2H_0,
\]
where $\nu_{\rm out}$ is the unit normal pointing
away from the center $x$.

\vspace{5pt}
\noindent Put $\gamma:=\Sigma'\cap\partial B_t(x)\cap U$,$D^\pm:=B_d^2(\pm de_3)\subset\mathbb R^2$ and
$\pi(p_1,p_2,p_3):=(p_2,p_3)$,then 
\begin{equation} \label{eq:surface intersect with sphere}
T_0(\gamma)=\operatorname{span}(e_2),\qquad \operatorname{spt}\|V''\|\cap\Gamma\cap U=\gamma.
\end{equation}
 where the second part follows, by virtue of \cite[Lemma~4.5]{DRKS20}, from the maximum principle in \cite[Proposition~2(iii)]{DDH19}. Choose $d>0$ sufficiently small so that
$\overline{D^\pm}\cap\pi(\gamma)=\{0\}$ and fix $\mu_0,\mu_1$ such that
\[
H_0<\mu_0<\kappa,
\qquad
\mu_1>
\sup_{\partial B_t(x)\cap U}
H_\Phi[\partial B_t(x),\nu_{\rm out}].
\]
Let $\mathcal M_\Phi$ denote the anisotropic
graph mean-curvature operator with the normal
directed towards negative $p_1$ direction.
By standard elliptic theory and the implicit
function theorem (cf. \cite[proof of
Proposition~6]{DDH19}), after shrinking
$d>0$ if necessary,there exist solutions $u_\mu^\pm$,
continuous in $\mu\in[\mu_0,\mu_1]$, of
\[
\begin{cases}
\mathcal M_\Phi(u_\mu^\pm)=\mu
&\text{in }D^\pm,\\
u_\mu^\pm=\psi
&\text{on }\partial D^\pm.
\end{cases}
\]
Then Hopf-maximum  principle yields
$u_{\mu_1}^\pm<\psi$ and $ u_{\mu_0}^\pm>\psi$ in 
$D^\pm.$As $\mu$ decreases from $\mu_1$ to $\mu_0$,
maximum principle rules out any first interior touches between
$\operatorname{graph}(u_\mu^\pm)$  and $\Sigma''$
. Consequently,
\begin{equation} \label{eq:barrier bound}
\operatorname{spt}\|V''\|
\cap B_t(x)\cap\pi^{-1}(D^\pm)
\subset
\{p_1\ge u_{\mu_0}^\pm(p_2,p_3)\}.
\end{equation}

\noindent By the Hopf boundary point lemma,there exist $a_\pm>0$ with
$\nabla(u_{\mu_0}^\pm-\psi)(0)=\pm a_\pm e_3$,denote $a=\min\{a_+,a_-\}$.
Let $r_j\downarrow0$ be such that
$(\eta_{y,r_j})_\#V''\rightharpoonup Z,$
since
\[
\begin{gathered}
r_j^{-1}\overline{D^\pm}
\longrightarrow\{\pm p_3\ge0\},\\
\frac{\psi(r_jq)}{r_j}\longrightarrow0,\qquad
\frac{u_{\mu_0}^\pm(r_jq)-\psi(r_jq)}{r_j}
\longrightarrow\pm a_\pm q_3,
\end{gathered}
\]
\eqref{eq:surface intersect with sphere} and \eqref{eq:barrier bound}
then implies
$\operatorname{spt}\|Z\|\cap\{p_1\ge0\}\subset\mathcal W_{1/a}$,
taking $L:=1/a$ proves the first part of \eqref{eq:pmc-spherical-wedge}.

Suppose $R_Z\neq0$ and set $K:=\operatorname{spt}\|R_Z\|\subset W_L,$
write $F=\Phi_y$ and
$F_{ij}(\xi)
:=\frac{\partial^2 F(\xi)}
{\partial \xi_i\partial \xi_j},
 1\le i,j\le3.$Fix $\kappa>0$ such that $\kappa L<1$.
For $\varepsilon\in(0,1)$, define
\[
f_\varepsilon(p)
=p_1-\kappa\sqrt{1+p_3^2}
+\varepsilon p_2^2.
\]
Then for every $p\in K$,
\[
\begin{aligned}
f_\varepsilon(p)
&\ge p_1-\kappa(1+|p_3|)
+\varepsilon p_2^2
\ge (1-\kappa L)p_1
+\varepsilon p_2^2-\kappa.
\end{aligned}
\]
Consequently, $f_\varepsilon(p)\to+\infty$ as
$|p|\to\infty$ with $p\in K$.
Hence, $f_\varepsilon$ attains its minimum on $K$
at some point $q_\varepsilon\in K$.
We denote this minimum by
$h_\varepsilon:=f_\varepsilon(q_\varepsilon)
=\min_{p\in K}f_\varepsilon(p)$.Set
\[
u_\varepsilon(p_2,p_3)
:=h_\varepsilon
+\kappa\sqrt{1+p_3^2}
-\varepsilon p_2^2.
\]
Then
\[
K\subset\{p_1\ge u_\varepsilon(p_2,p_3)\},
\qquad
(q_\varepsilon)_1
=u_\varepsilon((q_\varepsilon)_2,
               (q_\varepsilon)_3),
\]
and
\[
D^2u_\varepsilon=
\begin{pmatrix}
-2\varepsilon&0\\
0&\displaystyle
\frac{\kappa}{(1+p_3^2)^{3/2}}
\end{pmatrix}.
\]

Fix $q_0\in K$. From
$f_\varepsilon(q_\varepsilon)
\le f_\varepsilon(q_0)\le f_1(q_0)$
and the preceding lower bound for $f_\epsilon$, there exists
$C>0$, independent of $\varepsilon$, such that
\begin{equation}\label{eq:wedge-contact-bounds}
\begin{gathered}
0\le(q_\varepsilon)_1\le C,\\
|(q_\varepsilon)_3|\le LC,\qquad
\varepsilon(q_\varepsilon)_2^2\le C.
\end{gathered}
\end{equation}
Set
$N_\varepsilon
:=\{p\in\mathbb R^3:
p_1\ge u_\varepsilon(p_2,p_3)\},
S_\varepsilon:=\partial N_\varepsilon.$
Then
$K\subset N_\varepsilon,
q_\varepsilon\in K\cap S_\varepsilon.$Let $\nu_\varepsilon$ be the unit normal to $S_\varepsilon$ pointing into $N_\varepsilon$,namely
\[
\nu_\varepsilon
=\frac{(1,-Du_\varepsilon)}
{\sqrt{1+|Du_\varepsilon|^2}}
\]
By \eqref{eq:wedge-contact-bounds},
\[
|Du_\varepsilon(q_\varepsilon)|^2
\le\kappa^2+4C\varepsilon,
\]
and hence
\[
1-\langle\nu_\varepsilon,e_3\rangle^2
\ge
\frac{1}{1+\kappa^2+4C}.
\]

By uniform ellipticity of $F$,
\[
F_{33}(\nu_\varepsilon)\ge\lambda^{-1}\bigl(
1-\langle\nu_\varepsilon,e_3\rangle^2
\bigr)
\ge\frac{\lambda^{-1}}{1+\kappa^2+4C},
\qquad
F_{22}(\nu_\varepsilon)\le\lambda.
\]
Thus, by \cite[Equation~(11)]{DDH19},
\[
\begin{aligned}
&F(\nu_\varepsilon)
H_F(S_\varepsilon,\nu_\varepsilon)|_{q_\varepsilon}
\\
&=
\frac{
-2\varepsilon F_{22}(\nu_\varepsilon)
+\dfrac{\kappa F_{33}(\nu_\varepsilon)}
{(1+(q_\varepsilon)_3^2)^{3/2}}
}
{\sqrt{1+|Du_\varepsilon(q_\varepsilon)|^2}}\\
&\ge
\frac{1}{\sqrt{1+|Du_\varepsilon(q_\varepsilon)|^2}}
\left[
\frac{\kappa\lambda^{-1}}
{(1+L^2C^2)^{3/2}(1+\kappa^2+4C)}
-2\lambda\varepsilon
\right]\\
&>0
\end{aligned}
\]
for all sufficiently small $\varepsilon>0$.

On the other hand, since $\delta_F R_Z=0$,
\cite[Lemma~4.5]{DRKS20} implies that
$K=\operatorname{spt}\|R_Z\|$ satisfies
the maximum principle condition required in
\cite[Proposition~2(iii)]{DDH19}.
Applying the latter at $q_\varepsilon$, we obtain
\[
F(\nu_\varepsilon)
H_F(S_\varepsilon,\nu_\varepsilon)|_{q_\varepsilon}\le 0
\]
a contradiction. Hence $R_Z=0$.

\end{proof}

Thus we have proved 
$\operatorname{TanVar}(y,V'')=\{|T_y\Sigma'|\}.$
Apply the same graphicality argument as in \cite[Proof of Proposition 4.14,Step 1,Case2]{DD24} yields
\[
\Sigma''\cap B_t(x)\cap B_{r_y}(y)
=\Sigma'\cap B_t(x)\cap B_{r_y}(y)
\]
for some $r_y>0$, completing the case when $y\in \mathcal{R}(\Sigma')$.

\medskip
\noindent\textit{Case 2: $y\in \mathcal{S}(\Sigma')$.}
Since by \eqref{eq:count touching point} $\# \mathcal{S}(\Sigma')\cap\partial B_t(x)$ is finite,we can choose $r>0$
such that $B_{2r}(y)\Subset\operatorname{An}(x;\rho,2\rho)$ and
\[
\mathcal{S}(\Sigma')\cap\partial B_t(x)\cap B_{2r}(y)=\{y\},
\qquad
\Sigma'\cap B_{2r}(y)=\Sigma'_+\cup\Sigma'_-.
\]
Here $\Sigma'_\pm$ are smooth embedded sheets with common tangent
plane $P:=T_y\Sigma'$ at $y$. Put $H_0:=\|h\|_{C^0(M)}$ and
$\Phi_y:=\Phi(y,\cdot)$. Choose coordinates with
$T_y\partial B_t(x)=\{p_1=0\}$, and write
$H^\pm:=\{\pm p_1>0\}$.
The mass-ratio bound and the fact
$V''=|\Sigma'|$ outside $B_t(x)$ yield
\[
\operatorname{TanVar}(y,V'')\neq\varnothing,
\qquad
Z\llcorner H^-=2|P|\llcorner H^-
\quad\text{for every }Z\in\operatorname{TanVar}(y,V'').
\]

Fix such a $Z$ and $r_j\downarrow 0$ with
$W_j:=(\eta_{y,r_j})_\#V''\rightharpoonup Z.$
For $0<\alpha<1$, set
$A_\alpha:=
B_{1/\alpha}(0)\setminus\overline{B_\alpha(0)}.$
By the conclusion of Case 1 and unique continuation
along $\Sigma'$ we have
$V''\llcorner(B_{2r}(y)\setminus\{y\})\geq
|\Sigma'|\llcorner(B_{2r}(y)\setminus\{y\}).$
Consequently, for every fixed $\alpha\in(0,1)$ and
all sufficiently large $j$,
$W_j\llcorner A_\alpha
\geq
(\eta_{y,r_j})_\#|\Sigma'|
\llcorner A_\alpha.$
Since
$(\eta_{y,r_j})_\#|\Sigma'|
\rightharpoonup 2|P|,$passing to the limit yields
$Z\llcorner A_\alpha
\geq 2|P|\llcorner A_\alpha.$
As $\alpha$ is arbitrary, we obtain
$Z\llcorner(\mathbb R^3\setminus\{0\})
\geq
2|P|\llcorner(\mathbb R^3\setminus\{0\}).$
Moreover, since
$\|Z\|(B_a(0))\leq Ca^2$ for all sufficiently
small $a>0$, we have $\|Z\|(\{0\})=0$,
thus we conclude that $Z\geq 2|P|$.

Finally, since both $Z$ and $2|P|$ are stationary
for the frozen integrand $F=\Phi(y,\cdot)$ and 
$Z\llcorner H^-=2|P|\llcorner H^-,$ we have
\begin{equation}\label{eq:pmc-touching-remainder}
R_Z:=Z-2|P|\geq0,
\qquad
\delta_F R_Z=0,
\qquad
\operatorname{spt}\|R_Z\|
\subset\overline{H^+}.
\end{equation}

Apply the same argument in the proof of Lemma~\ref{lem:pmc-interface-wedge}
 gives $R_Z=0$. Thus
\begin{equation}\label{eq:pmc-touching-tangent}
\operatorname{TanVar}(y,V'')=\{2|T_y\Sigma'|\}.
\end{equation}

Now by the same graphicality argument as in \cite[Proof of Proposition 4.14,Step 1,Case2]{DD24} yields
\[
\Sigma''\cap B_t(x)\cap B_{r_y}(y)
=\Sigma'\cap B_t(x)\cap B_{r_y}(y)
\]
for some $r_y>0$, completing the case when $y\in \mathcal{S}(\Sigma')$.
This completes Step~1.

\medskip
\noindent\textbf{Step 2.}
We prove that there is a smooth almost embedded $h$-stable
APMC boundary $\Sigma\subset B_t(x)\setminus\{x\}$ such that
\begin{equation}\label{eq:pmc-step2-identification}
V\llcorner(B_\rho(x)\setminus\{x\})
=|\Sigma|\llcorner(B_\rho(x)\setminus\{x\}).
\end{equation}
Choose the initial radius $\rho$ sufficiently small that
$H_\Phi[\partial B_a(x),\nu_{\rm out}]>2H_0$ for
$0<a\leq2\rho$, and that 
\cite[Lemma~5.6]{DD24} applies with $c=H_0$.

For $0<s<\rho$, let $W_s$ be an $h$-replacement of $V'$ in
$\operatorname{An}(x;s,t)$,where $t$ as in Step 1.Let $\Sigma_s:=W_s|_{An(x;s,t)}$,then we have
\begin{equation}\label{eq:boundary-intersection}
\varnothing\neq\overline C\cap\partial B_t(x)
\subset\Sigma'\cap\partial B_t(x)
\quad\text{for every connected component }C\text{ of }\Sigma_s.
\end{equation}

Indeed, a component $C$  failing to reach $\partial B_t(x)$ has
an maximum of $d_g(x,\cdot)$ in $C$.Consequently, there exists some
$r_3<t$
such that $C$ touches $\partial B_{r_3}(x)$ from the inside. More
precisely,for some sufficiently small $\varepsilon>0$,
\[
C\cap\partial B_{r_3}(x)\neq\emptyset,\quad C\cap \bigl[\overline{B_{r_3+\varepsilon}(x)}\setminus\overline{B_{r_3}(x)}\bigr]=\emptyset.
\]
But, by \cite[Lemma~4.5]{DRKS20}, this contradicts \cite[Proposition~2(iii)]{DDH19}.

We observe that Step 1, combined with unique continuation and \eqref{eq:boundary-intersection}, yields

\[
\Sigma_{s_1}\cap\operatorname{An}(x;s_2,t)=\Sigma_{s_2}
\quad(0<s_1<s_2<\rho)
\]
Now define $\Sigma:=\bigcup_{0<s<\rho}\Sigma_s,$
thus $\Sigma$ is smooth, almost embedded stable $h-$APMC boundary
in $B_t(x)\setminus\{x\}$,
the boundary structure and $h$-stability of $\Sigma$
follow from those of the compatible local replacements.
Moreover,
\begin{equation}\label{eq:pmc-step2-replacement-sides}
W_s\llcorner B_s(x)=V\llcorner B_s(x),
\qquad
W_s\llcorner\operatorname{An}(x;s,t)
=|\Sigma|\llcorner\operatorname{An}(x;s,t).
\end{equation}

Since $V$ is an integral $2$-varifold by \cite[Lemma~4.15]{DD24}, \cite[Lemma~5.6]{DD24} implies that the set
\[
\mathcal D:=\left\{
q\in\operatorname{spt}\|V\|\cap\big(B_\rho(x)\setminus\{x\}\big):
\begin{array}{l}
\operatorname{TanVar}(q,V)=\{m_q|P_q|\},\quad m_q>0,\\
P_q\text{ is transverse to }T_q\partial B_{d_g(x,q)}(x)
\end{array}\right\}
\]
is dense in $\operatorname{spt}\|V\|\cap\big(B_\rho(x)\setminus\{x\}\big)$.
We remark that the cited lemma requires only $H_0$-bounded first variation,
and therefore applies in the present setting.

For $q\in\mathcal D$, put $a=d_g(x,q)$. Transversality and
\eqref{eq:pmc-step2-replacement-sides} imply
\[
q\in\spt\|V\llcorner B_a(x)\|
\subset\spt\|W_a\|.
\]
If $q\notin\overline{\spt\|W_a\|\setminus\overline{B_a(x)}}$,
then $\spt\|W_a\|$ locally touches $\partial B_a(x)$ from
inside, by \cite[Lemma~4.5]{DRKS20}, this contradicts \cite[Proposition~2(iii)]{DDH19}. Therefore
$q\in\overline{\Sigma\cap\operatorname{An}(x;a,t)}\subset\Sigma$.
Density of $\mathcal{D}$ yields
\begin{equation}\label{eq:pmc-step2-support}
\spt\|V\|\cap(B_\rho(x)\setminus\{x\})\subset\Sigma.
\end{equation}

It remains to determine the multiplicities.To this end,we establish the following constancy theorem.We remark that the classical constancy theorem for stationary varifolds
\cite[Theorem~41.1]{Sim83} and its anisotropic counterpart
\cite[Proposition~5]{DDH19} assume vanishing first variation and
the following variant allows nonzero anisotropic first
variation, provided the varifold is integral and it has bounded anisotropic first variation.
\begin{proposition}

\label{prop:bounded-constancy}
Let $U\subset M$ be open and $\Phi$ be an elliptic integrand in $U$, and let $\mathcal S^k\subset U$ be a connected, $C^2$, properly embedded submanifold without boundary in $U$.
 Let $W\in\mathcal{IV}_{k}(U)$ satisfy
$\spt\|W\|\subset S$ and has $c$ bounded $\Phi$-first variation in $U$

Then
\[
W=m|S|
\qquad\text{for some }m\in\mathbb N\cup\{0\}.
\]
\end{proposition}
\begin{proof}
The support inclusion give
$W=\theta|S|$, with
$\theta\in L^1_{\mathrm{loc}}(S;\mathbb N\cup\{0\})$; the
approximate tangent plane is $T_xS$ for $\|W\|$-almost every $x$.
Write $a(x)=\Phi(x,T_xS)>0$.  If $X$ is a compactly supported
$C^1$ tangent vector field on $S$, extend it to $U$.  Its flow
preserves $S$, thus
\begin{equation}
 \delta_\Phi W(X)=\int_S\theta\,
          \operatorname{div}_S(aX)\,\dd\mathcal H_g^2.
 \label{eq:tangential-constancy}
\end{equation}

On a relatively compact chart of $S$, set $X=Y/a$ for an arbitrary
compactly supported tangent field $Y$.  The variation bound implies
\[
 \left|\int_S\theta\operatorname{div}_S Y\,\dd\mathcal H_g^2\right|
 \leq C\int_S\theta|Y|_g\,\dd\mathcal H_g^2.
\]
Consequently the distributional derivative $D_S\theta$ is a Radon
measure with $|D_S\theta|\leq C\theta\mathcal H_g^2$,which implies $\theta\in W^{1,1}_{\rm loc}(S)$.
A Sobolev function with integer values has zero weak gradient,to see this,
apply the  chain rule to
$F(t)=(2\pi)^{-1}\sin(2\pi t)$; since $F(\theta)=0$ and
$F'(\theta)=1$ almost everywhere, $\nabla\theta=0$.
Thus $\theta\in \mathbb{N}\cup \{0\}$ is constant by connectedness of $S$. 
\end{proof}

We now return to the proof of Theorem~\ref{thm:APMC-regularity}
and determine the multiplicities.
Let $p\in \mathcal{R}(\Sigma)$ and $U\Subset B_\rho(x)\setminus\{x\}$ be an open neighborhood
of  $p$, chosen so that
$S:=\Sigma\cap U$ is connected and embedded.
By \eqref{eq:pmc-step2-support} and
Proposition~\ref{prop:bounded-constancy},
$V\llcorner U=m|S|$ for some $m\in\mathbb N\cup\{0\}$.
Choose a point $q\in S$ transverse to
$\partial B_a(x)$, $a=d_g(x,q)$.

Set
$U_-:=U\cap B_a(x)$ and
$U_+:=U\setminus\overline{B_a(x)}$.
Then \eqref{eq:pmc-step2-replacement-sides} gives
\[
W_a\llcorner U_-=m|S\cap U_-|,
\qquad W_a\llcorner U_+=|S\cap U_+|.
\]
Any support point of $W_a$ in
$(U\cap\partial B_a(x))\setminus S$ is again excluded by
\cite[Proposition~2(iii)]{DDH19}. Consequently
\[
\spt\|W_a\|\cap U\subset S,
\qquad
|\delta_\Phi(W_a\llcorner U)(X)|
\leq H_0\int_U|X|\,\dd\|W_a\|
\quad(X\in C_c^1(U;TM)).
\]
Proposition~\ref{prop:bounded-constancy} now gives
$W_a\llcorner U=k|S|$. Thus we have $m=k=1$
and  $V=|\Sigma|$ near every $p\in \mathcal{R}(\Sigma)$.
Finally, $\mathcal H^2(\mathcal S(\Sigma)) = 0$ by Corollary~\ref{cor:apmc-boundary-touching}, and hence
$\|V\|(\mathcal S(\Sigma)) = \lVert |\Sigma| \rVert(\mathcal S(\Sigma)) = 0$.
Together with \eqref{eq:pmc-step2-support}, this establishes
\eqref{eq:pmc-step2-identification}.

\medskip
\noindent\textbf{Step 3.}
It remains to show $x$ is a removable singularity of $V$.Recall $V$ possesses a uniform mass ratio upper bound and by step 2 $V$ is a smooth almost embedded $h$-stable APMC boundary in $B_\rho (x) \setminus \{x\}$,then \eqref{eq:pmc-stability-prelim} of this paper for \cite[(2.5)]{DDL24}
,the same argument as \cite[Step~3 of Theorem~5.1]{DDL24} yields that  $x$ is a removable singularity.
\end{proof}

Having established the regularity of the almost-minimizing
varifolds, we are now ready to complete the proof of
Theorem~\ref{thm:apmc-minmax}.

\begin{proof}[Proof of Theorem~\ref{thm:apmc-minmax}]
By Lemma~\ref{lem:strict-endpoint-gap},
$\Wph>0$, which implies that $\Sigma$ is nontrivial.Then the theorem follows immediately by combining
Theorem~\ref{thm:AM-extraction-417} and
Theorem~\ref{thm:APMC-regularity}.

\end{proof}

Combining Theorem~\ref{thm:MR-ONVP} and
Theorem~\ref{thm:apmc-minmax}, we obtain the following corollary.

\begin{corollary}
\label{cor:APMC-minmax}
Let $(M^3,g)$ be a closed  Riemannian manifold,
let $\Phi$ be a smooth elliptic integrand on $M$ satisfying
\eqref{eq:elliptic-class}, and let $h\in C^\infty(M)$ be
$\Phi$-good with $\int_M h\,\dd V_g\geq 0$.
Then there exists $E\in\Cacc$ such that $\Sigma=\partial E$
is a nonempty, smooth, closed, almost embedded surface with
prescribed anisotropic mean curvature $h$ with respect to
the outward unit normal of $E$. Moreover,
\begin{equation}
\Aph(E)=\Wph,
\qquad
\operatorname{MR}_{\Phi,r_0}(|\Sigma|)\leq C,
\label{eq:APMC-output}
\end{equation}
where $C>0$ depends only on $(M,g)$, $\lambda$ and
$\|h\|_{C^0(M)}$, and
\[
\Theta^2(|\Sigma|,x)=
\begin{cases}
1, & x\in\mathcal{R}(\Sigma),\\
2, & x\in\mathcal{S}(\Sigma).
\end{cases}
\]
\end{corollary}

\begin{proof}[Proof of Theorem~\ref{thm:apmc-existence}]
The theorem follows from Proposition~\ref{thm:good-generic}
and Corollary~\ref{cor:APMC-minmax}, noting that every smooth
elliptic integrand satisfies \eqref{eq:elliptic-class} for
some $\lambda\geq 1$.
If $\int_M h\,\dd V_g<0$, we apply
Corollary~\ref{cor:APMC-minmax} to the $\Phi$-good function
$-h$ and take the complement of the resulting Caccioppoli set.
\end{proof}

\appendix
\section{Derivation of the anisotropic mean-curvature formula}
\label{app:first-variation-formula}

We derive \eqref{eq:Hformula} with the sign convention in
\eqref{eq:H-convention}, using the horizontal--vertical derivative and
covariant-derivative notation fixed in
\cref{subsec:anisotropic-pmc}.  If $x(t)$ is a smooth curve and
$v(t)\in T_{x(t)}M\setminus\{0\}$ is a smooth vector field along it,
the connection splitting of $T(TM)$ gives the chain rule
\begin{equation}
\left.\frac{\dd}{\dd t}\right|_{t=0}\tPhi(x(t),v(t))
=d_x^{\mathrm h}\tPhi|_{(x(0),v(0))}(\dot x(0))
+d_v\tPhi|_{(x(0),v(0))}
\left(\left.\frac{D^g v}{\dd t}\right|_{t=0}\right).
\label{eq:horizontal-vertical-chain-rule}
\end{equation}

Let $\Sigma\subset M$ be a smooth closed two-sided surface with unit
normal $\nu$.  Let
$X$ be a smooth vector field on $M$, let $F_t$ be its flow, and set
$\Sigma_t=F_t(\Sigma)$.  Thus
\[
F_0=\operatorname{Id}_M,
\qquad
\left.\frac{\partial F_t}{\partial t}\right|_{t=0}=X.
\]
For $x\in\Sigma$, write $x_t=F_t(x)$, and let
$\widehat\nu_t(x)=\nu_t(x_t)\in T_{x_t}M$, where $\nu_t$ is the
smoothly chosen unit normal to $\Sigma_t$ with $\nu_0=\nu$. Put
\[
\dot\nu(x)=\left.\frac{D^g}{\dd t}\right|_{t=0}
\widehat\nu_t(x)\in T_xM.
\]
If $e_1,e_2$ are an orthonormal basis of $T_x\Sigma$, choose a curve
$\gamma(s)\subset\Sigma$ with $\gamma(0)=x$ and
$\gamma'(0)=e_\alpha$, and consider
$\mathcal F(s,t)=F_t(\gamma(s))$. Writing $D_s,D_t$ for
covariant differentiation along this variation, torsion freeness
and the fact $\partial_t\mathcal F=X\circ\mathcal F$ give
\begin{equation}
\left.\frac{D^g E_\alpha}{dt}\right|_{t=0}
=
\bigl(D_t \partial_s \mathcal{F}\bigr)(0,0)
=
\bigl(D_s \partial_t \mathcal{F}\bigr)(0,0)
=
\nabla^g_{e_\alpha} X.
\label{eq:transported-tangent-derivative}
\end{equation}

\noindent Since $\widehat\nu_t(x)$ is a unit normal to $\Sigma_t$, we have
\[
g_{x_t}(\widehat\nu_t(x),E_\alpha(t))=0,
\qquad
g_{x_t}(\widehat\nu_t(x),\widehat\nu_t(x))=1.
\]

\noindent Differentiating at $t=0$ and using \eqref{eq:transported-tangent-derivative} yields
\begin{equation}
g_x(\dot\nu,e_\alpha)
=-g_x(\nu,\nabla_{e_\alpha}X),
\qquad
g_x(\dot\nu,\nu)=0.
\label{eq:material-normal-derivative}
\end{equation}

Let $J_\Sigma F_t$ denote the tangential Jacobian of $F_t|_\Sigma$.
For
\[
G_{\alpha\beta}(t)
=g_{x_t}((dF_t)_x(e_\alpha),(dF_t)_x(e_\beta)),
\]
one has $J_\Sigma F_t=(\det G(t))^{1/2}$ and
\[
\dot G_{\alpha\beta}(0)
=g_x(\nabla_{e_\alpha}X,e_\beta)
+g_x(e_\alpha,\nabla_{e_\beta}X).
\]
Consequently,
\begin{equation}
\left.\frac{\dd}{\dd t}\right|_{t=0}J_\Sigma F_t
=\sum_{\alpha=1}^2g_x(\nabla_{e_\alpha}X,e_\alpha)
=\operatorname{div}_\Sigma X.
\label{eq:tangential-Jacobian-derivative}
\end{equation}
Since $\tPhi(x,\nu)=\Phi(x,\nu)$ for $|\nu|_g=1$, the change-of-variables
formula, \eqref{eq:horizontal-vertical-chain-rule}, and
\eqref{eq:tangential-Jacobian-derivative} yield
\begin{align}
\delta_\Phi|\Sigma|(X)
=\int_\Sigma\bigl(&d_x^{\mathrm h}\tPhi|_{(x,\nu)}(X)
+d_v\tPhi|_{(x,\nu)}(\dot\nu)
+\Phi(x,\nu)\operatorname{div}_\Sigma X\bigr)
\,\dd\mathcal H_g^2.
\label{eq:unintegrated-anisotropic-first-variation}
\end{align}

Now take $X=f\nu$ along $\Sigma$, using any smooth ambient extension.
Equation \eqref{eq:material-normal-derivative} gives
$\dot\nu=-\nabla_\Sigma f$, while
$\operatorname{div}_\Sigma(f\nu)=
f\operatorname{div}_\Sigma\nu$. We adopt the following notational convention.
\[
\xi_\Phi(x)=D_v\tPhi(x,\nu(x))\in T_xM,
\qquad
\xi_\Phi^\top=\xi_\Phi-g_x(\xi_\Phi,\nu)\nu
\in T_x\Sigma.
\]
Then \eqref{eq:unintegrated-anisotropic-first-variation} becomes
\begin{align}
\delta_\Phi|\Sigma|(f\nu)
=\int_\Sigma\bigl(&
f\,d_x^{\mathrm h}\tPhi|_{(x,\nu)}(\nu)
-g_x(\xi_\Phi^\top,\nabla_\Sigma f)
+f\Phi(x,\nu)\operatorname{div}_\Sigma\nu\bigr)
\,\dd\mathcal H_g^2.
\label{eq:normal-first-variation-before-integration}
\end{align}
Because $\Sigma$ is closed, tangential integration by parts gives
\[
-\int_\Sigma g_x(\xi_\Phi^\top,\nabla_\Sigma f)
\,\dd\mathcal H_g^2
=\int_\Sigma f\operatorname{div}_\Sigma(\xi_\Phi^\top)
\,\dd\mathcal H_g^2.
\]
Positive one-homogeneity of $\tPhi$ and Euler's identity give
$g_x(\xi_\Phi,\nu)=\Phi(x,\nu)$.  Hence
$\xi_\Phi=\xi_\Phi^\top+\Phi(x,\nu)\nu$, and therefore
\begin{equation}
\operatorname{div}_\Sigma\xi_\Phi
=\operatorname{div}_\Sigma(\xi_\Phi^\top)
+\Phi(x,\nu)\operatorname{div}_\Sigma\nu.
\label{eq:Cahn-Hoffman-divergence-splitting}
\end{equation}

Combining the preceding identities gives
\[
\delta_\Phi|\Sigma|(f\nu)
=\int_\Sigma f\left(
\operatorname{div}_\Sigma(D_v\tPhi(x,\nu))
+d_x^{\mathrm h}\tPhi|_{(x,\nu)}(\nu)
\right)\,\dd\mathcal H_g^2.
\]
Finally,
since $d_x^{\mathrm h}\tPhi|_{(x,\nu)}(\nu)
=g_x(D_x\tPhi(x,\nu),\nu)$,we have

\[
H_\Phi[\Sigma,\nu]
=\operatorname{div}_\Sigma(D_v\tPhi(x,\nu))
+g_x(D_x\tPhi(x,\nu),\nu),
\]
which is \eqref{eq:Hformula}.
\section{Allard-type regularity under a uniform lower mass-ratio bound}
\label{app:allard-uniform-mass-lower}

We record the following Allard-type regularity theorem, referring to
Allard's regularity result in \cite{All86}.
The hypotheses below are not formulated exactly as in Allard's original
statement; however, as indicated in the remarks in \cite{All86},
the theorem remains valid in the present form.

\begin{theorem}[Allard-type regularity under a uniform lower mass-ratio bound]
\label{thm:allard-uniform-mass-lower}

Let $n\ge 2$, $\alpha\in(0,1)$, $\lambda>1$, and
$M,\kappa,s>0$. Let $\Omega\subset\mathbb R^{n+1}$
be a fixed open set containing $\{0\}\times[-s,s]$,
and let $\Phi$ be a smooth elliptic parametric
integrand on $\Omega$ satisfying
\eqref{eq:elliptic-class} with constant $\lambda$.

Writing $z=(x,t)\in\mathbb R^n\times\mathbb R$, set
\[
\mathcal C_\rho
:=\Omega\cap\bigl(B_\rho^n(0)\times\mathbb R\bigr),
\qquad
\omega_n:=\mathcal L^n(B_1^n(0)).
\]

There exist positive constants $r_*(n,\alpha,M,\lambda),\delta_0(n,\alpha,M,\kappa,\lambda)$, and $C(n,\alpha,M,\kappa,\lambda)$,
,
with the following property.

Suppose $0<r<r_*$ and $V\in\mathbf{IV}_n(\Omega)$
satisfy the following conditions:

\begin{enumerate}

\item[(i)] 
$B_{4r}^n(0)\times[-s,s]\subset\Omega$, and
$
\operatorname{spt}\|V\|\cap\mathcal C_{4r}
\subset\{(x,t):|t|\le s\}.
$

\item[(ii)] 
$(4r)^{-n}\|V\|(\mathcal C_{4r})\le M$.

\item[(iii)] 
\begin{equation}\label{eq:allard-mass-interval}
\frac{\omega_n}{2}
<
(2r)^{-n}\|V\|(\mathcal C_{2r})
\le\frac{3\omega_n}{2}.
\end{equation}

\item[(iv)] 
There exists $H\ge0$ such that
$\|\delta_\Phi V\|\le H\|V\|$ on $\mathcal C_{4r}$.

\item[(v)] 
For every
$z\in\operatorname{spt}\|V\|\cap\mathcal C_{3r}$
and $0<\rho\le r/4$, 
\begin{equation}\label{eq:allard-mass-lower}
\|V\|(B_\rho^{n+1}(z)\cap\Omega)
\ge\kappa\rho^n.
\end{equation}

\item[(vi)] 
\begin{equation}\label{eq:allard-smallness}
\mu+4rH\le\delta_0.
\end{equation}
where
\[
\mu:=
\left(
(4r)^{-n-2}
\int_{\mathcal C_{4r}}
t^2\,d\|V\|(x,t)
\right)^{1/2},
\]

\end{enumerate}

\par\medskip\noindent
\emph{Then} there exists
$f\in C^{1,\alpha}(B_{2r}^n(0))$ such that
\begin{equation}\label{eq:allard-graph}
\operatorname{spt}\|V\|\cap\mathcal C_{2r}
=
\{(x,f(x)):x\in B_{2r}^n(0)\}.
\end{equation}

Moreover, with $\mathcal E:=\mu+rH+r$, we have
\begin{equation}\label{eq:allard-center-estimate}
\frac{|f(0)|}{r}+|Df(0)|
\le C\mathcal E
\end{equation}
and
\begin{equation}\label{eq:allard-holder-estimate}
|Df(\widetilde x)-Df(x)|
\le C\mathcal E
\left(\frac{|\widetilde x-x|}{r}\right)^\alpha
\end{equation}
for all $x,\widetilde x\in B_{2r}^n(0)$.

\end{theorem}


\medskip
\noindent\textbf{Declaration of Generative AI Use.}
The authors used ChatGPT during the preparation of this manuscript to assist
with refining the language, clarity, and presentation of the mathematical
exposition. The authors reviewed and revised
the resulting content as necessary and take full responsibility for the final
version of the article.

\vspace{2em} 

\noindent
Xiaoxiang Jiao\\
School of Mathematical Sciences, University of Chinese Academy of Sciences\\
19A Yuquan Road, Beijing, 100049, China\\
E-mail: xxjiao@ucas.ac.cn

\vspace{1.5em} 

\noindent
Wenduo Zou\\
School of Mathematical Sciences, University of Chinese Academy of Sciences\\
19A Yuquan Road, Beijing, 100049, China\\
E-mail: zouwenduo23@mails.ucas.ac.cn
\end{document}